\documentclass[11pt]{article}
\usepackage{amsthm, amsmath, amssymb, amsfonts, url, booktabs, tikz, setspace, fancyhdr, amsbsy}
\usepackage[margin = 0.8in]{geometry}
\usepackage{esint}
\usepackage{verbatim}

\definecolor{refkey}{gray}{.75}
\definecolor{labelkey}{gray}{.2}

\newtheorem{theorem}{Theorem}[section]
\newtheorem{proposition}[theorem]{Proposition}
\newtheorem{lemma}[theorem]{Lemma}

\theoremstyle{definition}

\theoremstyle{remark}
\newtheorem*{remark}{Remark}

\newcommand{\norm}[1]{\left\lVert#1\right\rVert}

\newcommand{\R}{\mathbb{R}}
\newcommand{\V}{\mathbb{V}}
\newcommand{\Q}{\mathbb{Q}}
\newcommand{\Z}{\mathbb{Z}}
\newcommand{\defeq}{\mathrel{\mathop:}=}
\newcommand{\eqdef}{\mathrel{\mathop=}:}

\newcommand{\G}{\mathcal{G}}
\newcommand{\HH}{\mathcal{H}}
\newcommand{\uu}{\boldsymbol{u}}
\newcommand{\DD}{\boldsymbol{D}}
\newcommand{\UU}{\boldsymbol{U}}
\newcommand{\VV}{\boldsymbol{V}}
\newcommand{\vv}{\boldsymbol{v}}
\newcommand{\SSS}{\boldsymbol{S}}
\newcommand{\q}{\boldsymbol{q}}
\newcommand{\ww}{\boldsymbol{w}}
\newcommand{\hh}{\boldsymbol{h}}
\newcommand{\dx}{\,\mathrm{d}x}

\numberwithin{equation}{section}

\def\ocirc#1{\ifmmode\setbox0=\hbox{$#1$}\dimen0=\ht0
    \advance\dimen0 by1pt\rlap{\hbox to\wd0{\hss\raise\dimen0
    \hbox{\hskip.2em$\scriptscriptstyle\circ$}\hss}}#1\else
    {\accent"17 #1}\fi}

\begin{document}

\title{Finite element approximation of steady flows of generalized
              Newtonian fluids with concentration-dependent power-law index}

\author{Seungchan Ko\thanks{Mathematical Institute, University of Oxford, Andrew Wiles Building, Woodstock Road, Oxford OX2 6GG, UK. Email: \tt{seungchan.ko@maths.ox.ac.uk}}
~and ~Endre S\"uli\thanks{Mathematical Institute, University of Oxford, Andrew Wiles Building, Woodstock Road, Oxford OX2 6GG, UK. Email: \tt{endre.suli@maths.ox.ac.uk}}}

\date{~}

\maketitle

~\vspace{-1.5cm}

\begin{abstract}
We consider a system of nonlinear partial differential equations describing the motion of an incompressible chemically reacting
generalized Newtonian fluid in three space dimensions. The governing system consists of a steady convection-diffusion equation
for the concentration and a  generalized steady power-law-type fluid flow model for the velocity and the pressure, where the
viscosity depends on both the shear-rate and the concentration through a concentration-dependent power-law index.
The aim of the paper is to perform a mathematical analysis of a finite element approximation of this model. We formulate a
regularization of the model by introducing an additional term in the conservation-of-momentum equation and construct a
finite element approximation of the regularized system. We show the convergence of the finite element method to a weak solution
of the regularized model and prove that weak solutions of the regularized problem converge to a weak solution of the
original problem.
\end{abstract}


\smallskip


\begin{section}{Introduction}
We are interested in developing a convergence theory for finite element approximations of a system of nonlinear partial differential
equations (PDEs) modelling the rheological response of the synovial fluid. The synovial fluid is a biological fluid found in the
cavities of movable joints and is composed of ultrafiltrated blood, called {\it{hyaluronan}}. Laboratory experiments have shown that
the viscosity of the fluid depends on the concentration of hyaluronan, as well as on the shear-rate. In particular, it was observed
in steady shear experiments that the concentration of the hyaluronan is not just a scaling factor of the viscosity (understood as
$\nu(c,|Du|)=f(c)\,\tilde{\nu}(|Du|)$) but it has an influence on the degree of shear-thinning. Therefore, a new mathematical model
of the rheological response of the synovial fluid was proposed in \cite{exp2}. There, the authors considered a power-law-type model
for the velocity and the pressure, where the power-law index depends on the concentration, corresponding to the fact that the
concentration affects the level of shear-thinning. To close the system, a generalized convection-diffusion equation was assumed to
be satisfied by the concentration. For a detailed rheological background we refer to \cite{exp2, exp}.

Based on the description above, we consider the following system of PDEs:
\begin{alignat}{2}
{\rm{div}}\,\uu&=0\qquad &&{\rm{in}}\,\,\Omega,\label{eq1}\\
{\rm{div}}\,(\uu\otimes\uu)-{\rm{div}}\,\SSS(c,\DD\uu)&=-\nabla p+\boldsymbol{f}\qquad &&{\rm{in}}\,\,\Omega,\label{eq2}\\
{\rm{div}}\,(c\uu)-{\rm{div}}\,\q_c(c,\nabla c, \DD\uu)&=0\qquad &&{\rm{in}}\,\,\Omega,\label{eq3}
\end{alignat}
where $\Omega\subset\R^d$ is a bounded open Lipschitz domain. In the above system of PDEs,
$\uu:\overline{\Omega}\rightarrow\R^d$, $p:\Omega\rightarrow\R$, $c:\overline{\Omega}\rightarrow\R_{\geq0}$ denote the velocity,
pressure and concentration fields, respectively, $\boldsymbol{f}:\Omega\rightarrow\R^d$ is a given external force, and $\DD\uu$ denotes
the symmetric velocity gradient, i.e. $\DD\uu=\frac{1}{2}(\nabla\uu+(\nabla\uu)^T)$. To complete the model, we impose the following
Dirichlet boundary conditions:
\begin{equation}\label{bc}
\uu=\mathbf{0},\qquad c=c_d\qquad{\rm{on}}\,\,\partial\Omega,
\end{equation}
where $c_d\in W^{1,s}(\Omega)$ for some $s>d$. By Sobolev embedding, $c_d$ is continuous up to the boundary, and we can therefore define
\[c^-\defeq\min_{x\in\overline{\Omega}}c_d\,\,\,\,\,\text{and}\,\,\,\,\,c^+\defeq\max_{x\in\overline{\Omega}}c_d.\]

We further assume that the extra stress tensor $\SSS:\R_{\geq 0}\times\R^{d\times
d}_{\rm sym}\rightarrow\R^{d\times d}_{\rm sym}$ is a continuous mapping
satisfying the following growth, strict monotonicity and coercivity
conditions, respectively: there exist positive constants $C_1$, $C_2$ and $C_3$ such that
\begin{equation}\label{S1}
|\SSS(\xi,\boldsymbol{B})|\leq C_1(|\boldsymbol{B}|^{r(\xi)-1}+1),
\end{equation}
\begin{equation}\label{S2}
(\SSS(\xi,\boldsymbol{B_1})-\SSS(\xi,\boldsymbol{B_2}))\cdot(\boldsymbol{B_1}-\boldsymbol{B_2})>0\,\,\,\text{for}\,\,\boldsymbol{B_1}\neq
\boldsymbol{B_2},
\end{equation}
\begin{equation}\label{S3}
\SSS(\xi,\boldsymbol{B})\cdot \boldsymbol{B}\geq C_2(|\boldsymbol{B}|^{r(\xi)}+|\SSS|^{r'(\xi)})-C_3,
\end{equation}
where $r:\R_{\geq0}\rightarrow\R_{\geq0}$ is a H\"older-continuous function satisfying $1<r^-\leq r(\xi)\leq r^+<\infty$ and $r'(\xi)$
is defined as its H\"older conjugate, $\frac{r(\xi)}{r(\xi)-1}$. We further assume that the concentration flux vector
$\q_c(\xi,\boldsymbol{g},\boldsymbol{B}):\R_{\geq 0}\times\R^d\times\R^{d\times d}_{\rm{sym}}\rightarrow\R^d$ is a
continuous mapping, which is linear with respect to $\boldsymbol{g}$, and it additionally satisfies the following growth and coercivity
conditions: there exist positive constants $C_4$ and $C_5$ such that
\begin{align}
|\q_c(\xi,\boldsymbol{g},\boldsymbol{B})|&\leq C_4|\boldsymbol{g}|,\label{q1}\\
\q_c(\xi,\boldsymbol{g},\boldsymbol{B})\cdot \boldsymbol{g}&\geq C_5|\boldsymbol{g}|^2.\label{q2}
\end{align}

As we have discussed above, the prototypical examples we have in mind are the following:
\[\SSS(c,\DD\uu)=\nu(c,|\DD\uu|)\DD\uu,\qquad \q_c(c,\nabla c,\DD\uu)=\boldsymbol{K}(c,|\DD\uu|)\nabla c,\]
where the viscosity $\nu(c,|\DD\uu|)$ is of the following form:
\[\nu(c,|\DD\uu|)\sim\nu_0(\kappa_1+\kappa_2|\DD\uu|^2)^{\frac{r(c)-2}{2}}\]
and $\nu_0,\kappa_1,\kappa_2$ are positive constants.

The rigorous mathematical analysis of the existence of global weak solutions to a PDE system, consisting of the generalized
Navier--Stokes equations, with a concentration-dependent viscosity coefficient, coupled to a convection-diffusion equation,
was initiated in \cite{BMR2008}. There, however, the power-law index was fixed and the concentration was only a scaling factor of the
viscosity; the authors considered the evolutionary model and established the long-time existence of large-data global weak solutions.
Concerning the model \eqref{eq1}--\eqref{q2} where the power-law index is concentration-dependent, the mathematical analysis was
initiated in \cite{BP2013}. The authors established there the existence of weak solutions, provided that $r^->\frac{3d}{d+2}$,
by using generalized monotone operator theory. In \cite{BP2014}, with the help of a Lipschitz-truncation technique, the existence
of weak solutions with $r^->\frac{d}{2}$ was proved and the H\"older continuity of the concentration was shown by using De Giorgi's method.

In \cite{KS2017}, the convergence of a finite element approximation to the system \eqref{eq1}--\eqref{q2} was shown, using a discrete
De Giorgi regularity result. Because of the absence of a discrete De Giorgi regularity result in three space dimensions, the analysis in
\cite{KS2017} was restricted to the case of two space dimensions. In this paper, we extend the analysis developed in \cite{KS2017} to
three space dimensions, and we formulate an analogous convergence result for a finite element method in a three-dimensional domain.
The main idea here is to use a different numerical approximation scheme from the one in \cite{KS2017}, resulting in a different
limiting process. To this end, we consider different meshes for the conservation of linear momentum equation and the concentration
equation. The resulting numerical method can be viewed as a two-level Galerkin approximation. This enables us to
separate the passages to the limits with respect to the discretization parameters in the two equations, thus avoiding the need for
a discrete De Giorgi regularity result in three space dimensions.

As a first step, in Section 2 we introduce the necessary notational conventions and auxiliary results, which will be used throughout the paper.
In Section 3, we define a regularized problem, which enables us to enlarge the range of the power-law index so as to be able to cover the
practically relevant range of values of this index. In Sections 4 and 5, we construct a two-level Galerkin finite element approximation to
the regularized problem and perform a convergence analysis of the numerical method. Finally, in Section 6, we shall prove that weak
solutions of the regularized problem converge to a weak solution of the original problem when we pass to the limit with the regularization
parameter.
\end{section}
\begin{section}{Notation and auxiliary results}
In this section, we shall introduce certain function spaces and auxiliary results that will be used throughout the paper. Let $\mathcal{P}$ be
the set of all measurable functions $r:\Omega\rightarrow[1,\infty]$; we shall call the function $r\in\mathcal{P}(\Omega)$ a variable
exponent. We define {\color{black}{$r^-\defeq$ ess $\inf_{x\in\Omega}r(x)$, $r^+\defeq$ ess $\sup_{x\in\Omega}r(x)$}} and we only consider the case
\begin{equation}\label{pcondition}
1<r^-\leq r^+<\infty.
\end{equation}

Since we are considering a power-law index depending on the concentration, we need to work with Lebesgue and Sobolev spaces with
variable exponents. To be specific, we introduce the following variable-exponent Lebesgue spaces, equipped with the corresponding
Luxembourg norms:
\begin{align*}
L^{r(\cdot)}(\Omega)&\defeq \left\{u\in L^1_{\rm{loc}}(\Omega):\int_{\Omega}|u(x)|^{r(x)}\dx<\infty\right\},\\
\norm{u}_{L^{r(\cdot)}(\Omega)}=\norm{u}_{r(\cdot)}&\defeq \inf\left\{\lambda>0:\int_{\Omega}\bigg|\frac{u(x)}{\lambda}\bigg|^{r(x)}\dx
\leq1\right\}.
\end{align*}
Similarly, we introduce the following generalized Sobolev spaces:
\begin{align*}
W^{1,r(\cdot)}(\Omega)&\defeq \left\{u\in W^{1,1}(\Omega)\cap L^{r(\cdot)}(\Omega):|\nabla u|\in L^{r(\cdot)}\right\},\\
\norm{u}_{W^{1,r(\cdot)}(\Omega)}=\norm{u}_{1,r(\cdot)}&\defeq
\inf\left\{\lambda>0:\int_{\Omega}\left[\bigg|\frac{u(x)}{\lambda}\bigg|^{r(x)}+\bigg|\frac{\nabla u(x)}{\lambda}\bigg|^{r(x)}\right]\dx
\leq1\right\}.
\end{align*}
It is easy to show that all of the above spaces are Banach spaces, and because of \eqref{pcondition}, they are all separable and
reflexive; see \cite{DHHR2011}.

Furthermore, we introduce certain function spaces that are frequently used in PDE models of incompressible fluids. Henceforth,
$X(\Omega)^d$ will denote the space of $d$-component vector-valued functions with components from $X(\Omega)$. We also define
the space of tensor-valued functions $X(\Omega)^{d\times d}$. Finally, we define the following spaces:
\begin{align*}
W^{1,r(\cdot)}_0(\Omega)^d&\defeq \left\{\uu\in W^{1,r(\cdot)}(\Omega)^d:\uu=\mathbf{0}\,\,\text{on}\,\,\partial\Omega\right\},\\
W^{1,r(\cdot)}_{0,\rm{div}}(\Omega)^d&\defeq \left\{\uu\in W^{1,r(\cdot)}_0(\Omega)^d:{\rm{div}}\,\uu=0 \,\,\text{in $\Omega$}\right\},\\
L^{r(\cdot)}_0(\Omega)&\defeq \left\{f\in L^{r(\cdot)}(\Omega):\int_{\Omega}f(x)\dx=0\right\}.
\end{align*}

Throughout the paper, we shall denote the duality pairing between $f\in X$ and $g\in X^*$ by  $\langle g,f\rangle$, and for two
vectors $\boldsymbol{a}$ and $\boldsymbol{b}$, $\boldsymbol{a}\cdot \boldsymbol{b}$ denotes their scalar product; similarly, for
two tensors $\mathbf{A}$ and $\mathbf{B}$, $\mathbf{A} \cdot \mathbf{B}$ signifies their scalar product. Also, for any Lebesgue-measurable set $Q\subset\R^d$, $|Q|$ denotes the standard Lebesgue measure of the set $Q$.

Next we introduce the necessary technical tools. First we define the subset $\mathcal{P}^{\rm{log}}(\Omega)\subset\mathcal{P}(\Omega)$:
it will denote the set of all log-H\"older-continuous functions
defined on $\Omega$, that is the set of all functions $r \in \mathcal{P}(\Omega)$ satisfying
\begin{equation}\label{logh}
|r(x)-r(y)|\leq\frac{C_{\rm{log}}(r)}{-\log|x-y|}\qquad \forall\, x,y\in\Omega:0<|x-y|\leq\frac{1}{2}.
\end{equation}
It is obvious that classical H\"older-continuous functions on $\Omega$ automatically belong to this class.

Next we state the following lemma, which summarizes some inequalities involving variable-exponent norms.
{\color{black}{For proofs, see \cite{DHHR2011}, which is an extensive source of information concerning variable-exponent spaces.}}

\begin{lemma}
Let $\Omega\subset\R^d$ be a bounded open Lipschitz domain and let $r\in\mathcal{P}^{\rm{log}}(\Omega)$ satisfy \eqref{pcondition}.
Then, the following inequalities hold:
\begin{itemize}
\item H\"older's inequality, i.e.,
\[\norm{fg}_{s(\cdot)}\leq2\norm{f}_{r(\cdot)}\norm{g}_{q(\cdot)},\,\,\,\text{with}\,\,\,r,q,s\in\mathcal{P}(\Omega),
\,\,\,\frac{1}{s(x)}=\frac{1}{r(x)}+\frac{1}{q(x)}, \quad x \in \Omega.\]

\item Poincar\'e's inequality, i.e.,
\[\|u\|_{r(\cdot)}\leq C(d,C_{\rm{log}}(r))\,{\rm{diam}}(\Omega)\|\nabla u\|_{r(\cdot)}\qquad\forall\, u\in W^{1,r(\cdot)}_0(\Omega).\]

\item Korn's inequality, i.e.,
\[\norm{\nabla \uu}_{r(\cdot)}\leq C(\Omega, C_{\rm{log}}(r))\norm{\DD\uu}_{r(\cdot)}\qquad \forall\, \uu\in W^{1,r(\cdot)}_0(\Omega)^d,\]
where $C_{\rm{log}}(r)$ is the constant appearing in the definition of the class of log-H\"older-continuous functions.
\end{itemize}
\end{lemma}

Another important auxiliary result is the existence of the Bogovski\u{\i} operator in the variable-exponent setting.
\begin{theorem}\label{contibog}
Let $\Omega\subset\R^d$ be a bounded open Lipschitz domain and suppose that $r\in\mathcal{P}^{\rm{log}}(\Omega)$ with $1<r^-\leq r^+<\infty$.
Then, there exists a bounded linear operator $\mathcal{B}:L^{r(\cdot)}_0(\Omega)\rightarrow W^{1,r(\cdot)}_0(\Omega)^d$ such that for
all $f\in L^{r(\cdot)}_0(\Omega)$ we have
\begin{align*}
{\rm{div}}\,(\mathcal{B}f)&=f,\\
\|\mathcal{B}f\|_{1,r(\cdot)}&\leq C\|f\|_{r(\cdot)},
\end{align*}
where $C$ depends on $\Omega$, $r^-$, $r^+$, and $C_{\rm{log}}(r)$.
\end{theorem}

Let us now state the inf-sup condition, which has a crucial role in the mathematical analysis of incompressible fluid flow problems.
\begin{proposition}
For any $s$, $s'\in(1,\infty)$, with $\frac{1}{s}+\frac{1}{s'}=1$, there exists a positive constant $\alpha_s>0$ such that
\begin{equation}\label{contiinfsup}
\alpha_s\|q\|_{s'}\leq\sup_{0\neq\vv\in W^{1,s}_0(\Omega)^d}\frac{\langle{\rm{div}}\,\vv,q\rangle}{\|\vv\|_{1,s}}\qquad\forall\,
q\in L^{s'}_0(\Omega).
\end{equation}
\end{proposition}
This is a direct consequence of the existence of the Bogovski\u{\i} operator in spaces with fixed exponent, which is a special
case of Theorem \ref{contibog};  see \cite{B1979, DRS2010} for additional details.

Furthermore, we can prove the following inf-sup condition in spaces with variable-exponent norms, which will play an important role
in the subsequent analysis.
\begin{proposition}\label{contiinfsup2}
{\color{black}{Let $\Omega\subset\R^d$ be a bounded open Lipschitz domain and let $r\in\mathcal{P}^{\rm{log}}(\Omega)$ with $1<r^-\leq r^+<\infty$.
Then, there exists a constant $\alpha_r>0$ such that
\[\alpha_r\|q\|_{r'(\cdot)}\leq\sup_{0\neq \vv\in W^{1,r(\cdot)}_0(\Omega)^d} \frac{\langle{\rm{div}}\,\vv,q\rangle}{\|\vv\|_{1,r(\cdot)}}
\qquad\forall\, q\in L^{r'(\cdot)}_0(\Omega).\]}}
\end{proposition}

Proposition \ref{contiinfsup2} is a direct consequence of Theorem \ref{contibog} and the norm-conjugate formula stated in the following lemma.

\begin{lemma}
Let $r\in\mathcal{P}^{\rm{log}}(\Omega)$ be a variable exponent with $1<r^-\leq r^+<\infty$; then we have
\[\frac{1}{2}\|f\|_{r(\cdot)}\leq\sup_{g\in L^{r'(\cdot)}(\Omega),\,\,\|g\|_{r'(\cdot)}\leq1}\int_{\Omega}|f||g|\dx,\]
for all measurable functions $f \in  L^{r(\cdot)}(\Omega)$.
\end{lemma}

Finally, we recall the following well-known result due to De Giorgi and Nash \cite{DG1957, N1958}; see also \cite{BF2002} for its
application to the system of partial differential equations considered in the present paper.

\begin{theorem}\label{DeGiorgi}
Let $\Omega\subset\R^d$ be a Lipschitz domain and let $s>d$ be fixed. Suppose that $\boldsymbol{K}\in L^{\infty}(\Omega)^{d\times d}$
is uniformly elliptic with ellipticity constant $\lambda>0$. Then, there exists an $\alpha\in(0,1)$ such that, for any
$\boldsymbol{f}\in L^s(\Omega)^d$, $g\in L^{\frac{ds}{d+s}}(\Omega)$ and any $c_d\in W^{1,s}(\Omega)$, there exists a unique
$c\in W^{1,2}(\Omega)$ such that $c-c_d\in W^{1,2}_0(\Omega)\cap C^{0,\alpha}(\overline{\Omega})$ and
\[\int_{\Omega}\boldsymbol{K}\nabla c\cdot\nabla\varphi\dx=\int_{\Omega}\boldsymbol{f}\cdot\nabla\varphi\dx+\int_{\Omega}g\varphi\dx
\qquad\forall\,\varphi\in W^{1,2}_0(\Omega);\]
furthermore, the following uniform bound holds:
\[\|c\|_{W^{1,2}\cap C^{0,\alpha}}\leq C\left(\Omega,\lambda,s,\|\boldsymbol{K}\|_{\infty},{\color{black}{\|\boldsymbol{f}\|_s}},\|g\|_{\frac{ds}{d+s}},
\|c_d\|_{1,s}\right).\]
\end{theorem}

Using these notations, the weak formulation of the problem \eqref{eq1}--\eqref{q2} is as follows.

{\bf{Problem (Q).}} For $\boldsymbol{f}\in (W^{1,r^-}_0(\Omega)^d)^*$, $c_d\in W^{1,s}(\Omega)$, $s>d$, and a H\"older-continuous
function $r$, with $1<r^-\leq r(c)\leq r^+<\infty$ for all $c\in[c^-,c^+]$, find $(c-c_d)\in W^{1,2}_0(\Omega)\cap
C^{0,\alpha}(\overline{\Omega})$, for some $\alpha \in (0,1)$, $\uu\in W^{1,r(c)}_0(\Omega)^d$, $p\in L^{r'(c)}_0(\Omega)$ such that
\begin{alignat*}{2}
\int_{\Omega}\SSS(c,\DD\uu)\cdot\nabla\boldsymbol{\psi}-(\uu\otimes \uu)\cdot\nabla\boldsymbol{\psi}\dx-\langle{\rm{div}}\,
\boldsymbol{\psi},p\rangle&=\langle \boldsymbol{f},\boldsymbol{\psi}\rangle\qquad&&\forall\,\boldsymbol{\psi}\in W^{1,\infty}_0(\Omega)^d,\\
\int_{\Omega}q\,{\rm{div}}\,\uu\dx&=0\qquad &&\forall\, q\in L^{r'(c)}_0(\Omega),\\
\int_{\Omega}\q_c(c,\nabla c,\DD\uu)\cdot\nabla\varphi-c\uu\cdot\nabla\varphi\dx&=0\qquad &&\forall\,\varphi\in W^{1,2}_0(\Omega).
\end{alignat*}

Thanks to Proposition \ref{contiinfsup2}, we can restate {\bf{Problem (Q)}} in the following (equivalent) divergence-free setting.

{\bf{Problem (P).}} For $\boldsymbol{f}\in (W^{1,r^-}_0(\Omega)^d)^*$, $c_d\in W^{1,s}(\Omega)$, $s>d$, and a H\"older-continuous
function $r$, with $1<r^-\leq r(c)\leq r^+<\infty$ for all $c\in[c^-,c^+]$, find $(c-c_d)\in C^{0,\alpha}(\overline{\Omega})
\cap W^{1,2}_0(\Omega)$, $\uu\in W^{1,r(c)}_{0,{\rm{div}}}(\Omega)^d$, such that
\begin{align*}
\int_{\Omega}\SSS(c,\DD\uu)\cdot\nabla\boldsymbol{\psi}-(\uu\otimes \uu)\cdot\nabla\boldsymbol{\psi}\dx&=\langle \boldsymbol{f},\boldsymbol{\psi}\rangle&&\forall\,\boldsymbol{\psi}\in W^{1,\infty}_{0,{\rm{div}}}(\Omega)^d,\\
\int_{\Omega}\q_c(c,\nabla c,\DD\uu)\cdot\nabla\varphi-c\uu\cdot\nabla\varphi\dx&=0&&\forall\,\varphi\in W^{1,2}_0(\Omega).
\end{align*}

From now on, for simplicity, we shall restrict ourselves to the case of $d=3$. Our results can be however easily extended to
the case of any $d \geq 2$. We note in passing that since no uniqueness result is currently known for weak solutions of the
problem under consideration, we can only prove that a subsequence of the sequence of discrete solutions converges to a weak
solution of the problem.
\end{section}
\begin{section}{Regularization of the problem}
Before constructing the approximation of {\color{black}{{\bf{problem (Q)}}}} we shall formulate a regularized problem; it will then be the
regularized problem that will be approximated by a finite element method. We shall show that the sequence of finite element
approximations converges to a weak solution of the regularized problem, and that solutions of the regularized problem, in turn,
converge to a weak solution of {\color{black}{{\bf{problem (Q)}}}}. The reason for proceeding in this way is that direct approximation of {\color{black}{{\bf{problem (Q)}}}},
which bypasses the use of the regularized problem, necessitates the imposition of an unnaturally strong condition on the variable exponent
$r$ in the convergence analysis of the finite element method; the procedure that we describe below does not suffer from this
shortcoming.

Motivated by \cite{reg},
we shall utilize the following regularized problem, involving the regularization parameter $k\in\mathbb{N}$.
We choose a sufficiently large $t>0$, such that $r^->\frac{3}{2}>\frac{t}{t-2}.$ Then we seek a weak solution
$(\uu,p,c)\defeq(\uu^k,p^k,c^k)$ to
\begin{align}
{\rm{div}}\,\uu&=0\qquad &&{\rm{in}}\,\,\Omega,\label{eeq1}\\
{\rm{div}}\,(\uu\otimes\uu)-{\rm{div}}\,\SSS(c,\DD\uu)+{\color{black}{\frac{1}{k}|\uu|^{t-2}\uu}}&=
-\nabla p+\boldsymbol{f}\qquad &&{\rm{in}}\,\,\Omega,\label{eeq2}\\
{\rm{div}}\,(c\uu)-{\rm{div}}\,\q_c(c,\nabla c, \DD\uu)&=0\qquad &&{\rm{in}}\,\,\Omega,\label{eeq3}
\end{align}
Therefore, we consider the following regularized weak formulation.

{\bf{Problem (Q*).}} For $\boldsymbol{f}\in (W^{1,r^-}_0(\Omega)^3)^*$, $c_d\in W^{1,s}(\Omega)$, $s>3$, and a H\"older-continuous
function $r$, with $1<r^-\leq r(c)\leq r^+<\infty$ for all $c\in[c^-,c^+]$, and $r^->\frac{3}{2}>\frac{t}{t-2}$, $t>2$, find
$(c-c_d)\defeq(c^k-c_d)\in W^{1,2}_0(\Omega)
\cap C^{0,\alpha}(\overline{\Omega})$, for some $\alpha \in (0,1)$, $\uu\defeq\uu^k\in W^{1,r(c)}_0(\Omega)^3$, $p\defeq
p^k\in L^{r'(c)}_0(\Omega)$ such that
\begin{align}
\int_{\Omega}\SSS(c,\DD\uu)\cdot\nabla\boldsymbol{\psi}-(\uu\otimes \uu)\cdot\nabla\boldsymbol{\psi}+
{\color{black}{\frac{1}{k}|\uu|^{t-2}\uu}}\cdot\boldsymbol{\psi}\dx-\langle{\rm{div}}\,\boldsymbol{\psi},p\rangle&=
\langle \boldsymbol{f},\boldsymbol{\psi}\rangle\quad&&\forall\,\boldsymbol{\psi}\in W^{1,\infty}_0(\Omega)^3,\label{RGal1}\\
\int_{\Omega}q\,{\rm{div}}\,\uu\dx&=0\quad &&\forall\, q\in L^{r'(c)}_0(\Omega),\label{RGal2}\\
\int_{\Omega}\q_c(c,\nabla c,\DD\uu)\cdot\nabla\varphi-c\uu\cdot\nabla\varphi\dx&=0\quad &&\forall\,\varphi\in W^{1,2}_0(\Omega).\label{RGal3}
\end{align}

Again, by using Proposition \ref{contiinfsup2}, we can restate {\bf{Problem (Q*)}} in the following (equivalent) divergence-free setting:

{\bf{Problem (P*).}} For $\boldsymbol{f}\in (W^{1,r^-}_0(\Omega)^3)^*$, $c_d\in W^{1,s}(\Omega)$, $s>3$, and H\"older-continuous function
$r$, with $1<r^-\leq r(c)\leq r^+<\infty$ for all $c\in[c^-,c^+]$,  and $r^->\frac{3}{2}>\frac{t}{t-2}$, $t>2$, find $(c-c_d)\defeq(c^k-c_d)\in C^{0,\alpha}(\overline{\Omega})\cap
W^{1,2}_0(\Omega)$, $\uu\defeq\uu^k\in W^{1,r(c)}_{0,{\rm{div}}}(\Omega)^3$, such that
\begin{align}
\int_{\Omega}\SSS(c,\DD\uu)\cdot\nabla\boldsymbol{\psi}-(\uu\otimes \uu)\cdot\nabla\boldsymbol{\psi}+{\color{black}{\frac{1}{k}|\uu|^{t-2}\uu}}\cdot\boldsymbol{\psi}\dx&=\langle \boldsymbol{f},\boldsymbol{\psi}\rangle&&\forall\,\boldsymbol{\psi}\in W^{1,\infty}_{0,{\rm{div}}}(\Omega)^3,\label{RRGal1}\\
\int_{\Omega}\q_c(c,\nabla c,\DD\uu)\cdot\nabla\varphi-c\uu\cdot\nabla\varphi\dx&=0&&\forall\,\varphi\in W^{1,2}_0(\Omega).\label{RRGal2}
\end{align}

We shall formulate the finite element approximation of the regularized problem {\bf{Problem (Q*)}} in a three-dimensional domain;
the convergence analysis of the method
is presented in Section 4 and Section 5. In Section 6, we will prove that a sequence of weak solution triples
$\{(\uu^k,p^k,c^k)\}_{k \geq 1}$
of the regularized problem converges to a weak solution triple ($\uu,p,c$) of {\bf{Problem (Q)}}.
The latter result is recorded in our next theorem.

\begin{theorem}\label{mainthm2}
Suppose that $\Omega\subset\R^3$ is a convex polyhedral domain and $c_d\in W^{1,s}(\Omega)$ for some $s>3$.
Let us further assume that $r:\R_{\geq0}\rightarrow\R_{\geq0}$ is a H\"older-continuous function with $r^->\frac{3}{2}>\frac{t}{t-2}$, $t>2$,
and suppose that $\boldsymbol{f}\in (W^{1,r^-}_0(\Omega)^3)^*$.
Let $(\uu^k,p^k,c^k)$ be a weak solution of the regularized problem \eqref{eeq1}--\eqref{eeq3}.
Then, as $k\rightarrow\infty$, (a subsequence, not indicated, of) the sequence $\{(\uu^k,p^k,c^k)\}_{k \geq 1}$
converges to $(\uu,p,c)$ in the following sense:
\begin{align*}
\uu^k & \rightharpoonup \uu && {\rm{weakly}}\,\,{\rm{in}}\,\,W^{1,r^-}_{0,{\rm{div}}}(\Omega)^3, \\
c^k & \rightharpoonup c && {\rm{weakly}}\,\,{\rm{in}}\,\,W^{1,2}(\Omega),\\
c^k & \rightarrow c && {\rm{strongly}}\,\,{\rm{in}}\,\,C^{0,\alpha}(\overline{\Omega})\qquad{\rm{for}}\,\,\,{\rm{some}}\,\,\,\alpha\in(0,1),\\
p^k & \rightharpoonup p && {\rm{weakly}}\,\,{\rm{in}}\,\,L^{j'}(\Omega)\qquad\forall\, j>\max\{r^+,2\}.
\end{align*}
Furthermore, $(\uu,p,c)$ is a weak solution of the problem {\bf{Problem (Q*)}} stated in \eqref{eq1}--\eqref{eq3}.
\end{theorem}

\end{section}

\begin{section}{Finite element approximation}
\begin{subsection}{Finite element spaces}
Let $\{\G_n\}$, $\{\HH_m\}$  be families of shape-regular partitions of $\overline{\Omega}$ such that the following properties hold:
\begin{itemize}
\item {\bf{Affine equivalence}}: For each element $E\in \G_n$ (or $E\in\HH_m$) , there exists an invertible affine mapping
\[\boldsymbol{F}_E:E\rightarrow\hat{E},\]
where $\hat{E}$ is the standard reference $3$-simplex in $\R^3$.
\item {\bf{Shape-regularity}}: For any element $E\in \G_n$ (or $E\in\HH_m$), {\color{black}{the ratio of ${\rm{diam}}\,E$ to
the radius of the inscribed ball is bounded below uniformly by a positive constant, with respect to all $\G_n$ (or $\HH_m$) and
$n\in\mathbb{N}$ (or $m \in \mathbb{N}$).}}
\end{itemize}

For given partitions $\G_n$ and $\HH_m$, the finite element spaces are defined by
\begin{align*}
\V^n&=\V(\G_n)\defeq \{\VV\in C(\overline{\Omega})^3:\VV_{|E}\circ \boldsymbol{F}^{-1}_E\in\hat{\mathbb{P}}_{\V},E\in \G_n\,\,\text{and}\,\,\VV_{|\partial\Omega}=\boldsymbol{0}\},\\
\Q^n&=\Q(\G_n)\defeq \{Q\in L^{\infty}(\Omega):Q_{|E}\circ \boldsymbol{F}^{-1}_E\in\hat{\mathbb{P}}_{\Q},E\in \G_n\},\\
\Z^m&=\Z(\HH_m)\defeq\{Z\in C(\overline{\Omega}):Z_{|E}\circ \boldsymbol{F}^{-1}_E\in\hat{\mathbb{P}}_{\Z},E\in \HH_m\,\,\text{and}\,\,Z_{|\partial\Omega}=0\},
\end{align*}
where $\hat{\mathbb{P}}_{\V}\subset W^{1,\infty}(\hat{E})^3$, $\hat{\mathbb{P}}_{\Q}\subset L^{\infty}(\hat{E})$ and $\hat{\mathbb{P}}_{\Z}\subset W^{1,\infty}(\hat{E})$ are finite-dimensional linear subspaces.

We assume that $\V^n$ and $\Z^m$ have finite and locally supported bases; for example,
for each $n\in\mathbb{N}$ and $m\in\mathbb{N}$, there exists an $N_n\in\mathbb{N}$ and an $N_m\in\mathbb{N}$ such that
\[\V^n=\text{span}\{\VV^n_1,\ldots,\VV^n_{N_n}\},\]
\[\Z^m=\text{span}\{Z^m_1,\ldots,Z^m_{N_m}\},\]
and for each basis function $\VV^n_i$, $Z^m_j$, we have that if there exists an $E\in \G_n$ (respectively, $\HH_m$), with $\VV^n_i\neq0$ (respectively, $Z^m_j\neq0$) on $E$, then
\[\text{supp}\,\VV^n_i\subset\bigcup\{E'\in \G_n:E'\cap E\neq\emptyset\}\eqdef S_E.\]
\[\text{supp}\,Z^m_j\subset\bigcup\{E'\in \HH_m:E'\cap E\neq\emptyset\}\eqdef T_E.\]
For the pressure space $\Q^n$, we assume that $\Q^n$ has a basis consisting of discontinuous piecewise polynomials; i.e., for each $n\in\mathbb{N}$, there exists an $\tilde{N}_n\in\mathbb{N}$ such that
\[\Q^n={\rm{span}}\{Q^n_1,\ldots,Q^n_{\tilde{N}_n}\}\]
and for each basis function $Q^n_i$, we have that
\[{\rm{supp}}\,Q^n_i=E\qquad{\rm{for}}\,\,\,{\rm{some}}\,\,\,E\in\G_n.\]
We assume further that $\V^n$ contains continuous piecewise linear functions and $\Q^n$ contains piecewise constant functions.

Using the assumed shape-regularity we can easily verify that
\[\exists X\in\mathbb{N}:|S_E|\leq X|E|\,\,\,{\rm{for}}\,\,\,{\rm{all}}\,\,\,E\in\G_n,\]
\[\exists Y\in\mathbb{N}:{\color{black}{|T_E|}}\leq Y|E|\,\,\,{\rm{for}}\,\,\,{\rm{all}}\,\,\,E\in\HH_m,\]
where $X$ is independent of $n$ and $Y$ is independent of $m$. We denote by $g_E$ the diameter of $E\in\G_n$ and
by $h_E$ the diameter of $E\in\HH_m$.

We also introduce the subspace $\V^n_{\rm{div}}$ of discretely divergence-free functions. More precisely, we define
\[\V^n_{\rm{div}}\defeq\{\VV\in\V^n:\langle{\rm{div}}\,\VV,Q\rangle=0\,\,\,\forall\, Q\in\Q^n\},\]
and the subspace of $\Q^n$ consisting of vanishing integral mean-value approximations:
\[\Q^n_0\defeq\{Q\in\Q^n:\int_{\Omega}Q\dx=0\}.\]

Throughout this paper, we assume that the finite element spaces introduced above have the following minimal approximation properties.

\smallskip

\textbf{Assumption 1} (Approximability) For all $s\in[1,\infty)$,
\begin{align*}
\inf_{\VV\in\V^n}\norm{\vv-\VV}_{1,s}&\rightarrow0\qquad \qquad \forall\, \vv\in W^{1,s}_0(\Omega)^3\,\,\text{as}\,\,n\rightarrow\infty, \\
\inf_{Q\in\Q^n}\norm{q-Q}_s&\rightarrow0 \qquad \qquad \forall\, q\in L^s(\Omega)\,\,\text{as}\,\,n\rightarrow\infty,\\
\inf_{Z\in\Z^m}\norm{z-Z}_{1,s}&\rightarrow0\qquad \qquad \forall\, z\in W^{1,s}_0(\Omega)\,\,\text{as}\,\,m\rightarrow\infty.
\end{align*}
For this, a necessary condition is that the maximal mesh size vanishes, i.e., that $\max_{E\in\G_n}g_E\rightarrow0$ as $n\rightarrow\infty$ and $\max_{E\in\HH_m}h_E\rightarrow0$ as $m\rightarrow\infty$.

\smallskip

\textbf{Assumption 2} (Existence of a projection operator $\Pi^n_{\rm{div}}$) For each $n\in\mathbb{N}$, there exists a linear projection operator $\Pi^n_{\rm{div}}:W^{1,1}_0(\Omega)^3\rightarrow\V^n$ such that:
\begin{itemize}
\item $\Pi^n_{\rm{div}}$ preserves the divergence structure in the dual of the discrete pressure space; in other words, for any $\vv\in W^{1,1}_0(\Omega)^3$, we have
\[\langle{\rm{div}}\,\vv,Q\rangle=\langle{\rm{div}}\,\Pi^n_{\rm{div}}\vv,Q\rangle\qquad\forall\, Q\in\Q^n.\]
\item $\Pi^n_{\rm{div}}$ is locally $W^{1,1}$-stable, i.e., there exists a constant $c_1>0$, independent of $n$, such that
\begin{equation}\label{3.2}
\hspace{-4mm}\fint_E|\Pi^n_{\rm{div}}\vv|+g_E|\nabla\Pi^n_{\rm{div}}\vv|\dx\leq c_1\fint_{S_E}|\vv|+g_E|\nabla \vv|\dx\qquad\forall\, \vv\in W^{1,1}_0(\Omega)^3\,\,\,{\rm{and}}\,\,\,\forall\, E\in \G_n.
\end{equation}

\end{itemize}

Note that the local $W^{1,1}(\Omega)^3$-stability of $\Pi^n_{\rm{div}}$ implies its local and global $W^{1,s}(\Omega)^3$-stability for $s\in[1,\infty]$. In other words, for any $s\in[1,\infty]$ we have
\begin{equation}\label{sta_v}
\norm{\Pi^n_{\rm{div}} \vv}_{1,s}\leq c_s\norm{\vv}_{1,s}\qquad\forall\, \vv\in W^{1,s}_0(\Omega)^3,
\end{equation}
with a constant $c_s>0$ independent of $n>0$.

 Note further that the approximability ({\color{black}{{\bf{Assumption 1}}}}) and inequality \eqref{sta_v} imply the convergence of $\Pi^n_{\rm{div}} \vv$ to $\vv$. In fact, 
\begin{equation}\label{v_conv}
\|\vv-\Pi^n_{\rm{div}} \vv\|_{1,s}\rightarrow0\qquad\forall\, \vv\in W^{1,s}_0(\Omega)^3\,\,\text{as}\,\,n\rightarrow\infty, \quad
{\color{black}{\forall\, s \in [1,\infty)}}.
\end{equation}

\smallskip

\textbf{Assumption 3} (Existence of a projection operator $\Pi^n_{\Q}$) For each $n\in\mathbb{N}$,
there exists a linear projection operator $\Pi^n_{\Q}:L^1(\Omega)\rightarrow\Q^n$ such that $\Pi^n_{\Q}$ is locally $L^1$-stable; i.e.,
there exists a constant $c_2>0$, independent of $n$, such that
\begin{equation}\label{Q_L1sta}
\fint_E|\Pi^n_{\Q}q|\dx\leq c_2\fint_{S_E}|q|\dx
\end{equation}
for all $q\in L^1(\Omega)$ and all $E\in\G_n$.

\smallskip

Again, we have the following global stability and convergence property:
\begin{equation}\label{q_sta}
\|\Pi^n_{\Q}q\|_{s'}\leq c_{s'}\|q\|_{s'}\qquad \forall\, q\in L^{s'}(\Omega), \quad\forall\, s'\in(1,\infty),
\end{equation}
and
\begin{equation}\label{q_conv}
\|q-\Pi^n_{\Q}q\|_{s'}\rightarrow0,\qquad{\rm{as}}\,\,\,n\rightarrow\infty\,\,\,{\rm{for}}\,\,\,{\rm{all}}\,\,\,q\in L^{s'}(\Omega)\,\,\,{\rm{and}}\,\,\,s'\in(1,\infty).
\end{equation}

\begin{remark}
{\color{black}{According to \cite{BBDR2012}, the following pairs of velocity-pressure finite element spaces satisfy Assumptions 1, 2 and 3, for example:}}
\begin{itemize}
\item {\color{black}{The conforming Crouzeix--Raviart Stokes element, i.e., continuous piecewise quadratic plus cubic bubble velocity and discontinuous piecewise linear pressure approximation}} (compare e.g. with \cite{BF2012});
\item The space of continuous piecewise quadratic polynomials for the velocity and piecewise constant pressure approximation; see, \cite{BF2012}.
\end{itemize}
\end{remark}

Our final assumption is the existence of a projection operator for the concentration space.

\smallskip

\textbf{Assumption 4} (Existence of a projection operator {\color{black}{$\Pi^m_{\Z}$}}) For each $m\in\mathbb{N}$, there exists a linear projection operator $\Pi^m_{\Z}:W^{1,1}_0(\Omega)\rightarrow\Z^m$ such that
\[\fint_E|\Pi^m_{\Z}z|+h_E|\nabla\Pi^m_{\Z}z|\dx\leq c_3\fint_{T_E}|z|+h_E|\nabla z|\dx\qquad\forall\, z\in W^{1,1}_0(\Omega)\,\,\,{\rm{and}}\,\,\,\forall\, E\in \HH_m,\]
where $c_3$ does not depend on $m$.

Similarly as above, the projection operator $\Pi^m_{\Z}$ is globally $W^{1,s}$-stable for $s\in[1,\infty]$, and thus, by approximability,
\begin{equation}\label{z_conv}
\|\Pi^m_{\Z}z-z\|_{1,s}\rightarrow0 \qquad\forall\, z\in W^{1,s}_0(\Omega),\quad {\color{black}{\forall\, s \in [1,\infty)}}.
\end{equation}

Finally, we introduce a discrete inf-sup condition, which holds in our finite element setting.
It is a direct consequence of \eqref{contiinfsup} and the existence of $\Pi^n_{\rm{div}}$; see
\cite{BBDR2012} for further details.

\begin{proposition}\label{disinfsup1}
For $s$, $s'\in(1,\infty)$ satisfying $\frac{1}{s}+\frac{1}{s'}=1$, there exists a positive constant $\beta_r>0$, which is independent of $n$, such that
\[\beta_r\|Q\|_{s'}\leq\sup_{0\neq\VV\in \V^n}\frac{\langle{\rm{div}}\,\VV,Q\rangle}{\|\vv\|_{1,s}}\qquad\forall\, Q\in \Q^n_0\,\,\,{\rm{and}}\,\,\,\forall\,n\in\mathbb{N}.\]
\end{proposition}
\end{subsection}
\begin{subsection}{The finite element approximation}
In this section, we shall construct the finite element approximation of the problem \eqref{eeq1}--\eqref{eeq3}.
An important property of the incompressible Navier--Stokes equations is that the convective term in the momentum equation
is skew-symmetric; this is a consequence of the velocity field $\uu$ being divergence-free.
However, in the discretized problem, we might lose the skew-symmetry because we are considering only discretely divergence-free
finite element functions from the finite element space for the velocity.
Thus we need to modify the finite element approximation of the convective term in order to ensure that the skew-symmetry
is preserved under discretization. We therefore define the following modified convective terms:
\begin{align*}
B_u[\vv,\ww,\hh]&\defeq\frac{1}{2}\int_{\Omega}((\vv\otimes \hh)\cdot\nabla \ww-(\vv\otimes \ww)\cdot\nabla \hh)\dx,\\
B_c[b,\vv,z]&\defeq\frac{1}{2}\int_{\Omega}(z\vv\cdot\nabla b-b\vv\cdot\nabla z)\dx,
\end{align*}
for all $\vv,\ww,\hh\in W^{1,\infty}_0(\Omega)^3$, $b,z\in W^{1,\infty}(\Omega)$. These trilinear forms then
coincide with the corresponding trilinear forms appearing in the weak formulations of the momentum equation and the
concentration equation, provided that we are considering
pointwise divergence-free velocity fields. Furthermore, thanks to their skew symmetry, these two trilinear forms
now also vanish for discretely
divergence-free functions when $\ww=\hh$ and $b=z$, respectively. Explicitly, we have
\begin{align}
\begin{aligned}
B_u[\vv,\vv,\vv]&=0\,\,\,\,\,{\rm{and}}\,\,\,\,\,B_c[z,\vv,z]=0&&\forall\, \vv\in W^{1,\infty}_0(\Omega)^3,\,\,\,z\in W^{1,\infty}(\Omega),\label{Bfree}\\
B_u[\vv,\ww,\hh]&=-\int_{\Omega}(\vv\otimes \ww)\cdot\nabla \hh\dx&&\forall\, \vv,\ww,\hh\in W^{1,\infty}_{0,{\rm{div}}}(\Omega)^3,\\
B_c[b,\vv,z]&=-\int_{\Omega}b\vv\cdot\nabla z\dx&&\forall\, \vv\in W^{1,\infty}_{0,{\rm{div}}}(\Omega)^3,\,\,\,b,z\in W^{1,\infty}(\Omega).
\end{aligned}
\end{align}
Moreover, the trilinear form $B_u[\cdot,\cdot,\cdot]$ is bounded. Indeed, if $\vv,\ww,\hh\in W^{1,\infty}_0(\Omega)^3$,
then, by H\"older's inequality,
\[\int_{\Omega}(\vv\otimes\ww)\cdot\nabla\hh\dx\leq\|\vv\|_{2(r^-)'}\|\ww\|_{2(r^-)'}\|\hh\|_{1,r^-},\]
and
\[\int_{\Omega}(\vv\otimes\hh)\cdot\nabla\ww\dx\leq\|\vv\|_{2(r^-)'}\|\hh\|_{2(r^-)'}\|\ww\|_{1,r^-}.\]
Therefore, we obtain the bound
\begin{equation}\label{B_u1}
|B_u[\vv,\ww,\hh]|\leq \|\vv\|_{2(r^-)'}\|\ww\|_{2(r^-)'}\|\hh\|_{1,r^-}+\|\vv\|_{2(r^-)'}\|\ww\|_{1,r^-}\|\hh\|_{2(r^-)'}.
\end{equation}

Now, for each $n,m\in\mathbb{N}$, we call a triple $(\UU^{n,m},P^{n,m},C^{n,m})\in\V^n\times\Q^n_0\times(\Z^m+c_d)$ a discrete
solution to the Galerkin approximation if it satisfies
\begin{alignat}{2}
\int_{\Omega}\SSS(C^{n,m},\DD\UU^{n,m})\cdot\DD\VV+\frac{1}{k}|\UU^{n,m}|^{t-2}\UU^{n,m}\cdot\VV\dx&+B_u[\UU^{n,m},\UU^{n,m},\VV]\nonumber\\
-\langle{\rm{div}}\,\VV,P^{n,m}\rangle&=\langle\boldsymbol{f},\VV\rangle &&\forall\,\VV\in\V^n,\label{PGal1}\\
\int_{\Omega}Q\,{\rm{div}}\,\UU^{n,m}\dx&=0&&\forall\, Q\in\Q^n,\label{PGal2}\\
\int_{\Omega}\q_c(C^{n,m},\nabla C^{n,m},\DD\UU^{n,m})\cdot\nabla Z\dx+B_c[C^{n,m},\UU^{n,m},Z]&=0&&\forall\, Z\in\Z^m,\label{PGal3}
\end{alignat}
where $c_d\in W^{1,s}(\Omega)$ with $s>3$ and $\boldsymbol{f}\in(W^{1,r^-}_0(\Omega)^3)^*$.

If we restrict the test functions $\VV$ to $\V^n_{\rm{div}}$, then the above problem is transformed to the following:
find $(\UU^{n,m},C^{n,m})\in\V^n_{\rm{div}}\times(\Z^m+c_d)$ satisfying
\begin{align}
\int_{\Omega}\SSS(C^{n,m},\DD\UU^{n,m})\cdot\DD\VV+{\color{black}{\frac{1}{k}|\UU^{n,m}|^{t-2}}}\UU^{n,m}\cdot\VV\dx
+B_u[\UU^{n,m},\UU^{n,m},\VV]&=\langle\boldsymbol{f},\VV\rangle &&\forall\,\VV\in\V^n_{\rm{div}}\label{Gal1}\\
\int_{\Omega}\q_c(C^{n,m},\nabla C^{n,m},\DD\UU^{n,m})\cdot\nabla Z\dx+B_c[C^{n,m},\UU^{n,m},Z]&=0&&\forall\, Z\in\Z^m.\label{Gal2}
\end{align}

If $\frac{3}{2}<r^-$, the existence of the discrete solution pair $(\UU^{n,m}, C^{n,m})\in\V^n_{\rm{div}}\times(\Z^m+c_d)$
follows from a fixed point argument combined with an iteration scheme. Let us briefly summarize the proof of the existence
of the pair $(\UU^{n,m}, C^{n,m})\in\V^n_{\rm{div}}\times(\Z^m+c_d)$. Let $\{\ww_i\}^{N_n}_{i=1}$ be a basis of
$\V^n_{\rm{div}}\subset W^{1,\infty}_0(\Omega)^3$ such that $\int_{\Omega}\ww_i\cdot\ww_j\dx=\delta_{ij}$ and let
$\{z_j\}^{N_m}_{j=1}$ be a basis of $\Z^m\subset W^{1,2}_0(\Omega)$ such that $\int_{\Omega}z_iz_j=\delta_{ij}$.
Then, for fixed $n,m\in\mathbb{N}$, we define the Galerkin approximations.
\begin{equation}
\UU^{n,m}\defeq\sum^{N_n}_{i=1}\alpha_i^{n,m}\ww_i,\qquad C^{n,m}\defeq\sum^{N_m}_{i=1}\beta_i^{n,m}z_i+c_d,
\end{equation}
which satisfy \eqref{Gal1}--\eqref{Gal2}.

 First we define $C^{n,m}_1\defeq c_d\in\Z^m+c_d$. Then, for any $\ell\in\mathbb{N}$, we define $\UU^{n,m}_{\ell}\in\V^n_{\rm{div}}$
 as a solution of the finite-dimensional problem
\[\int_{\Omega}\SSS(C^{n,m}_{\ell},\DD\UU^{n,m}_{\ell})\cdot\DD\VV+{\color{black}{\frac{1}{k}|\UU^{n,m}_{\ell}|^{t-2}}}\UU^{n,m}\cdot\VV\dx+B_u[\UU^{n,m}_{\ell},\UU^{n,m}_{\ell},\VV]=\langle\boldsymbol{f},\VV\rangle \qquad\forall\,\VV\in\V^n_{\rm{div}},\]
and $C^{n,m}_{\ell}\in\Z^m+c_d$ as a solution of the finite-dimensional problem
\[\int_{\Omega}\q_c(C^{n,m}_{\ell},\nabla C^{n,m}_{\ell},\DD\UU^{n,m}_{\ell-1})\cdot\nabla Z\dx+B_c[C^{n,m}_{\ell},\UU^{n,m}_{\ell-1},Z]=0\qquad\forall\, Z\in\Z^m.\]
The existence of the functions  $\UU^{n,m}_{\ell}\in\V^n_{\rm{div}}$ and $C^{n,m}_{\ell}\in\Z^m+c_d$ is easily shown by means of Brouwer's fixed point theorem. Furthermore, for each $n,m\in\mathbb{N}$, the sequences of functions $\{\UU^{n,m}_{\ell}\}^{\infty}_{{\ell}=1}$ and $\{C^{n,m}_{\ell}\}^{\infty}_{{\ell}=1}$ satisfy the following uniform bounds:
\[\|\UU^{n,m}_{\ell}\|_{1,r^-}+{\color{black}{\|\UU^{n,m}_{\ell}\|_t}}\leq C_1,\qquad\|\nabla C^{n,m}_{\ell}\|_2\leq C_2,\]
where $C_1$ and $C_2$ are positive constants, independent of $\ell$.
Thus, by the Bolzano--Weierstrass theorem we deduce the existence of
limits $\UU^{n,m} \in \V^n_{\rm{div}}$ and $C^{n,m} \in \Z^m+c_d$ for $\UU^{n,m}_{\ell}$ and
$C^{n,m}_{\ell}$, respectively, as $\ell \rightarrow \infty$,
and these limits form a solution pair for the Galerkin approximation \eqref{Gal1}, \eqref{Gal2}.
For further details, see \cite{KS2017}.  This establishes the existence of a solution to the Galerkin approximations
\eqref{Gal1}, \eqref{Gal2} for any fixed pair of integers $n,m\in\mathbb{N}$.
The existence of a discrete solution triple for \eqref{PGal1}--\eqref{PGal3} then follows by the discrete inf-sup condition
stated in Proposition \ref{disinfsup1}, and we write $P^{n,m}=\sum^{\tilde{N}_n}_{i=1}\gamma^{n,m}_iy_i$
where $\{y_i\}^{\tilde{N}_n}_{i=1}$ is a basis of $\Q^n_0$.

We are now ready to state and prove our main theorem in this section. It asserts that, as $n,m\rightarrow\infty$, the sequence of discrete
solution triples converges to a weak solution triple of the regularized problem.
\begin{theorem}\label{mainthm1}
Suppose that $\Omega\subset\R^3$ is a convex polyhedral domain and $c_d\in W^{1,s}(\Omega)$ for some $s>3$.
Let us assume that $r:\R_{\geq0}\rightarrow\R_{\geq0}$ is a H\"older-continuous function with
$r^->\frac{3}{2}>\frac{t}{t-2}$, $t>2$, and let $\boldsymbol{f}\in (W^{1,r^-}_0(\Omega)^3)^*$.
Let $(\UU^{n,m},P^{n,m},C^{n,m})\in\V^n_{\rm{div}}\times\Q^n_0\times(\Z^m+c_d)$ be a discrete solution
triple defined by the finite element approximation \eqref{PGal1}--\eqref{PGal3}.
Then, the following convergence results hold.
\begin{itemize}
\item At the first level of Galerkin approximation, {\color{black}{there exists a subsequence (not relabelled) with respect to $m$ such that (as $m\rightarrow\infty$),}}
\begin{align*}
\UU^{n,m}&\rightarrow\UU^n&&{\mbox{uniformly on $\overline\Omega$}},\\
\DD\UU^{n,m}&\rightarrow\DD\UU^n&&{\mbox{uniformly on $\overline\Omega$}},\\
P^{n,m}&\rightarrow P^n&&{\mbox{uniformly on $\overline\Omega$}},\\
C^{n,m}&\rightharpoonup C^n&&{\rm{weakly}}\,\,{\rm{in}}\,\,W^{1,2}(\Omega),
\end{align*}
where $\UU^n\in\V^n$, $P^n\in\Q^n_0$.
\item At the second level of Galerkin approximation, {\color{black}{there exists a subsequence (not relabelled) with respect to $n$ such that (as $n\rightarrow\infty$),}}
\begin{align*}
\UU^n & \rightharpoonup \uu && {\rm{weakly}}\,\,{\rm{in}}\,\,W^{1,r^-}_{0}(\Omega)^3\\
P^n & \rightharpoonup p && {\rm{weakly}}\,\,{\rm{in}}\,\,L^{j'}(\Omega)\qquad\forall\, j>\max\{r^+,2\},\\
C^n & \rightharpoonup c && {\rm{weakly}}\,\,{\rm{in}}\,\,W^{1,2}(\Omega),\\
C^n & \rightarrow c && {\rm{strongly}}\,\,{\rm{in}}\,\,C^{0,\alpha}(\overline{\Omega}),\\
\end{align*}
where $(\uu,p,c)=(\uu^k,p^k,c^k)$ is a weak solution triple of the regularized problem \eqref{RGal1}--\eqref{RGal3}.
\end{itemize}
\end{theorem}
\end{subsection}
\end{section}

\begin{section}{Proof of Theorem \ref{mainthm1}}
\begin{subsection}{The limit $m\rightarrow\infty$}
First, we shall derive some uniform bounds, independent of $m\in\mathbb{N}$, and let $m$ tend
to infinity by using the weak compactness properties in the corresponding reflexive spaces. For simplicity, we shall denote $\SSS^{n,m}\defeq\SSS(C^{n,m},\DD\UU^{n,m})$, $\q_c^{n,m}\defeq\q_c(C^{n,m},\nabla C^{n,m},\DD\UU^{n,m})$.

We test with $\UU^{n,m}\in\V^{n}_{\rm{div}}$ in \eqref{PGal1}; then, thanks to the skew symmetry of $B_u[\cdot,\cdot,\cdot]$, we have
\[\int_{\Omega}\SSS^{n,m}\cdot\nabla\UU^{n,m}+\frac{1}{k}|\UU^{n,m}|^t\dx=\int_{\Omega}\SSS^{n,m}\cdot\DD\UU^{n,m}+\frac{1}{k}|\UU^{n,m}|^t\dx=\langle\boldsymbol{f},\UU^{n,m}\rangle.\]
By \eqref{S3} and Young's inequality, we have
\begin{equation}\label{uniform1}
\int_{\Omega}|\nabla\UU^{n,m}|^{r(C^{n,m})}+|\SSS^{n,m}|^{r'(C^{n,m})}+|\UU^{n,m}|^t\dx\leq C_1,
\end{equation}
where $C_1$ is independent of $m$.

Next, we test with $C^{n,m}-c_d \in \Z^m$ in \eqref{PGal3} and deduce that
\[\int_{\Omega}\q_c(C^{n,m},\nabla C^{n,m}\DD\UU^{n,m})\cdot\nabla(C^{n,m}-c_d)\dx=B_c[C^{n,m},\UU^{n,m},c_d].\]
By \eqref{q1}, \eqref{q2}, H\"older's inequality and Young's inequality,
\begin{align*}
\|\nabla C^{n,m}\|^2_2
&\leq\int_{\Omega}|\nabla C^{n,m}||\nabla c_d|\dx+B_c[C^{n,m},\UU^{n,m},c_d]\\
&\leq\varepsilon\|\nabla C^{n,m}\|^2_2+C(\varepsilon)\|\nabla c_d\|^2_2+B_c[C^{n,m},\UU^{n,m},c_d].
\end{align*}
Then, by Sobolev embedding,
\begin{align*}
B_c[C^{n,m},\UU^{n,m},c_d]
&=\frac{1}{2}\int_{\Omega}c_d\UU^{n,m}\cdot\nabla C^{n,m}\dx-\frac{1}{2}\int_{\Omega}C^{n,m}\UU^{n,m}\cdot\nabla c_d\dx\\
&=\int_{\Omega}c_d\UU^{n,m}\cdot\nabla C^{n,m}\dx+\frac{1}{2}\int_{\Omega}C^{n,m}({\rm{div}}\,\UU^{n,m})c_d\dx\\
&\leq\|c_d\|_{\infty}\|\UU^{n,m}\|_2\|\nabla C^{n,m}\|_2+\frac{1}{2}\,\|c_d\|_{\infty}\,\|C^{n,m}\|_{(r^-)'}\|{\rm{div}}\,\UU^{n,m}\|_{r^-}\\
&\leq C\|\UU^{n,m}\|_{1,r^-}\|\nabla C^{n,m}\|_2+C\|\UU^{n,m}\|_{1,r^-}\|\nabla C^{n,m}\|_{\frac{3r^-}{4r^--3}}\\
&\leq C(\varepsilon)\|\UU^{n,m}\|^2_{1,r^-}+\varepsilon\|\nabla C^{n,m}\|^2_2.
\end{align*}
Hence, by \eqref{q1} and \eqref{uniform1}, we have
\begin{equation}\label{uniform2}
\int_{\Omega}|\nabla C^{n,m}|^2+|\q^{n,m}_c|^2\dx\leq C_2,
\end{equation}
where $C_2$ is independent of $m$.

Next, we shall derive a uniform bound on the pressure. By Proposition \ref{disinfsup1} together with \eqref{PGal1}, \eqref{B_u1}
and the equivalence of norms in the finite-dimensional spaces, we have
\begin{align*}
\beta_r\|P^{n,m}\|_{(r^+)'}
&\leq\sup_{0\neq\VV\in\V^n}\frac{\langle{\rm{div}}\,\VV,P^{n,m}\rangle}{\|\VV\|_{1,r^+}}\\
&\leq\sup_{0\neq\VV\in\V^n}\frac{|\int_{\Omega}\SSS^{n,m}\cdot\DD\VV\dx|}{\|\VV\|_{1,r^+}}+C\sup_{0\neq\VV\in\V^n}\frac{|B_u[\UU^{n,m},\UU^{n,m},\VV]-\langle\boldsymbol{f},\VV\rangle|}{\|\VV\|_{1,r^-}}\\
&\leq C \sup_{0\neq\VV\in\V^n}\frac{\|\SSS^{n,m}\|_{(r^+)'}\|\DD\VV\|_{r^+}}{\|\VV\|_{1,r^+}}+C(n)\sup_{0\neq\VV\in\V^n}\frac{\|\UU^{n,m}\|^2_{2(r^-)'}\|\VV\|_{1,r^-}+\|\boldsymbol{f}\|_{-1}\|\VV\|_{1,r^-}}{\|\VV\|_{1,r^-}}.
\end{align*}
Therefore, by \eqref{uniform1}, we deduce that
\begin{equation}\label{uniform3}
\|P^{n,m}\|_{(r^+)'}\leq C(n).
\end{equation}

Now we are ready to let $m$ tend to infinity. By \eqref{uniform1} and \eqref{uniform3} with the equivalence of norms in
finite-dimensional spaces, we have $|\boldsymbol{\alpha}^{n,m}|\leq C(n)$ and $|\boldsymbol{\gamma}^{n,m}|\leq C(n)$.
Then, together with the uniform estimates \eqref{uniform2}, we can extract (not relabelled) subsequences such that
\begin{align}
\boldsymbol{\alpha}^{n,m}&\rightarrow\boldsymbol{\alpha}^n&&{\rm{strongly}}\,\,\,{\rm{in}}\,\,\,\R^{N_n},\label{conv1}\\
\boldsymbol{\gamma}^{n,m}&\rightarrow\boldsymbol{\gamma}^n&&{\rm{strongly}}\,\,\,{\rm{in}}\,\,\,\R^{\tilde{N}_n},\label{conv1-1}\\
C^{n,m}&\rightharpoonup C^n&&{\rm{weakly}}\,\,\,{\rm{in}}\,\,\,W^{1,2}(\Omega).\label{conv2}
\end{align}
From \eqref{conv1}, \eqref{conv1-1} and compact embedding,  we have
\begin{align}
\UU^{n,m}&\rightarrow\UU^n&&{\mbox{uniformly on $\overline\Omega$}},\label{conv6}\\
\DD\UU^{n,m}&\rightarrow\DD\UU^n&&{\mbox{uniformly on $\overline\Omega$}},\label{conv7}\\
P^{n,m}&\rightarrow P^n&&{\mbox{uniformly on $\overline\Omega$}},\label{convp}\\
C^{n,m}&\rightarrow C^n&&{\rm{strongly}}\,\,\,{\rm{in}}\,\,\,L^2(\Omega).\label{conv8}
\end{align}
By \eqref{conv1} and \eqref{conv1-1}, note that
\[\UU^n\in\V^n\qquad{\rm{and}}\qquad P^n\in\Q^n_0.\]
Finally, from \eqref{conv8}, we can extract a further subsequence (not relabelled) such that
\begin{equation}\label{conv9}
C^{n,m}\rightarrow C^n\qquad{\rm{a.e.}}\,\,\,{\rm{in}}\,\,\,\Omega.
\end{equation}

Note that since $\SSS$ is continuous, by \eqref{conv9} and \eqref{conv7}, we have
\[\SSS(C^{n,m},\DD\UU^{n,m})\rightarrow\SSS(C^n,\DD\UU^n)\qquad{\rm{a.e.}}\,\,\,{\rm{in}}\,\,\,\Omega.\]
Now, by \eqref{conv7}, we have that, for sufficiently large $m\in\mathbb{N}$,
\[|\DD\UU^{n,m}|<1+|\DD\UU^n|\qquad\mbox{for a.e. $x\in\Omega$}.\]
Thus, by \eqref{S1}, we have, for sufficiently large $m\in\mathbb{N}$,
\begin{align*}
|\SSS(C^{n,m},\DD\UU^{n,m})|
&\leq C|\DD\UU^{n,m}|^{r(C^{n,m})-1}+C\\
&\leq C(1+|\DD\UU^n|)^{r(c^{n,m})-1}+C\\
&\leq C(1+|\DD\UU^n|)^{r^+-1}+C,
\end{align*}
and $C(1+|\DD\UU^n|)^{r^+-1}+C\in L^{(r^+)'}(\Omega)$. Therefore, by the Dominated Convergence Theorem, we have
\begin{equation}\label{conv10}
\SSS^{n,m}\rightarrow\SSS^n\defeq\SSS(C^n,\DD\UU^n)\qquad{\rm{strongly}}\,\,\,{\rm{in}}\,\,\,L^{(r^+)'}(\Omega)^{3\times3}.
\end{equation}

Furthermore, by \eqref{conv9} and \eqref{conv7}, together with the Dominated Convergence Theorem,
\[\boldsymbol{K}(C^{n,m},|\DD\UU^{n,m}|)\rightarrow\boldsymbol{K}(C^n,|\DD\UU^n|)\qquad{\rm{strongly}}\,\,\,{\rm{in}}\,\,\,L^q(\Omega)\qquad\forall\, q\in(1,\infty).\]
Therefore, together with \eqref{conv2}, we have
\begin{equation}\label{conv11}
\q^{n,m}_c\rightharpoonup\q^n_c\defeq\q_c(C^n,\nabla C^n,\DD\UU^n)\qquad{\rm{weakly}}\,\,\,{\rm{in}}\,\,\,L^2(\Omega)^3.
\end{equation}

Now we are ready to pass $m$ to infinity in the Galerkin approximation \eqref{PGal1}--\eqref{PGal3}. First, by \eqref{conv6} and \eqref{conv7},
\begin{align*}
B_u[\UU^{n,m},\UU^{n,m},\VV]&\rightarrow B_u[\UU^n,\UU^n,\VV]&&\forall\, \VV\in\V^n,\\
\frac{1}{k}|\UU^{n,m}|^{t-2}\UU^{n,m}\cdot\DD\VV&\rightarrow\frac{1}{k}|\UU^n|^{t-2}\UU^n\cdot\DD\VV&&\forall\, \VV\in\V^n.
\end{align*}

Furthermore, from \eqref{conv10} and \eqref{convp},
\begin{align*}
\int_{\Omega}\SSS^{n,m}\cdot\DD\VV\dx&\rightarrow\int_{\Omega}\SSS^n\cdot\DD\VV\dx&&\forall\, \VV\in\V^n,\\
\langle{\rm{div}}\,\VV,P^{n,m}\rangle&\rightarrow\langle{\rm{div}}\,\VV,P^n\rangle&&\forall\, \VV\in\V^n.
\end{align*}
Therefore, we have
\begin{equation}\label{minfty1}
\int_{\Omega}\SSS^n\cdot\DD\VV+\frac{1}{k}|\UU^n|^{t-2}\UU^n\cdot\VV\dx+B_u[\UU^n,\UU^n,\VV]-\langle{\rm{div}}\,\VV,P^n\rangle=\langle\boldsymbol{f},\VV\rangle\qquad\forall \,\VV\in\V^n.
\end{equation}
Moreover, from \eqref{PGal2} and \eqref{conv7},
\begin{equation}\label{minfty2}
\int_{\Omega}Q\, {\rm{div}}\,\UU^n\dx=0\qquad\forall\, Q\in\Q^n.
\end{equation}

Next, let us investigate the limit of the concentration equation, \eqref{PGal3}. We fix an arbitrary $Z\in W^{1,2}_0(\Omega)$
and define $Z^m\defeq\Pi^m_{\Z}Z\in\Z^m$. Thanks to \eqref{conv6} and \eqref{conv8},
\begin{align*}
\|C^{n,m}\UU^{n,m}-C^n\UU^n\|_2
&\leq\|(\UU^{n,m}-\UU^n)C^{n,m}\|_2+\|\UU^n(C^{n,m}-C^n)\|_2\\
&\leq\|\UU^{n,m}-\UU^n\|_{\infty}\|C^{n,m}\|_2+\|\UU^n\|_{\infty}\|C^{n,m}-C^n\|_2\rightarrow0.
\end{align*}
Also, thanks to \eqref{conv6} and \eqref{z_conv},
\begin{align*}
\|Z^m\UU^{n,m}-Z\UU^n\|_2
&\leq\|(\UU^{n,m}-\UU^n)Z^m\|_2+\|\UU^n(Z^m-Z)\|_2\\
&\leq\|\UU^{n,m}-\UU^n\|_{\infty}\|Z^m\|_2+\|\UU^n\|_{\infty}\|Z^m-Z\|_2\rightarrow0.
\end{align*}
In other words, we have
\begin{align}
C^{n,m}\UU^{n,m}&\rightarrow C^n\UU^n&&{\rm{strongly}}\,\,\,{\rm{in}}\,\,\,L^2(\Omega)^3,\label{fconv1}\\
Z^m\UU^{n,m}&\rightarrow Z\UU^n&&{\rm{strongly}}\,\,\,{\rm{in}}\,\,\,L^2(\Omega)^3.\label{fconv2}
\end{align}
By \eqref{fconv2} and \eqref{conv2},
\begin{align*}
&\bigg|\int_{\Omega}Z^m\UU^{n,m}\cdot\nabla C^{n,m}\dx-\int_{\Omega}Z\UU^n\cdot\nabla C^n\dx\bigg|\\
&\quad\leq\int_{\Omega}|Z^m\UU^{n,m}-Z\UU^n||\nabla C^{n,m}|\dx+\bigg|\int_{\Omega}Z\UU^n(\nabla C^{n,m}-\nabla C^n)\dx\bigg|\rightarrow0.
\end{align*}
Moreover, from \eqref{fconv1} and \eqref{z_conv},
\begin{align*}
&\bigg|\int_{\Omega}C^{n,m}\UU^{n,m}\nabla Z^m\dx-\int_{\Omega}C^n\UU^n\cdot\nabla Z\dx\bigg|\\
&\quad\leq\|C^{n,m}\UU^{n,m}\|_2\|Z^m-Z\|_{1,2}+\|Z\|_{1,2}\|C^{n,m}\UU^{n,m}-C^n\UU^n\|_2\rightarrow0.
\end{align*}
Therefore, we have
\[\lim_{m\rightarrow\infty}B_c[C^{n,m},\UU^{n,m},Z^m]=B_c[C^n,\UU^n,Z].\]
Finally, from \eqref{conv11},
\[\int_{\Omega}\q^{n,m}_c\cdot\nabla Z^m\dx\rightarrow\int_{\Omega}\q^n_c\cdot\nabla Z\dx\qquad{\rm{as}}\,\,\,m\rightarrow\infty.\]
Altogether, we have
\begin{equation}\label{minfty3}
\int_{\Omega}\q^n_c\cdot\nabla Z\dx+B_c[C^n,\UU^n,Z]=0\qquad\forall\, Z\in W^{1,2}_0(\Omega).
\end{equation}
\end{subsection}

\begin{subsection}{The limit $n\rightarrow\infty$}
Now we shall derive further uniform estimates and let $n$ pass to infinity. First, we test with $\UU^n$ in \eqref{minfty1}.
Then, by \eqref{Bfree} and \eqref{minfty2}, we have
\[\int_{\Omega}\SSS^n\cdot\DD\UU^n+\frac{1}{k}|\UU^n|^t\dx=\langle\boldsymbol{f},\UU^n\rangle.\]
By using \eqref{S3} and Young's inequality, we have
\begin{equation}\label{uni1}
\int_{\Omega}|\DD\UU^n|^{r(C^n)}+|\SSS^n|^{r'(C^n)}+\frac{1}{k}|\UU^n|^t\dx\leq C_1,
\end{equation}
where $C_1$ is independent of $n$, which leads us to
\begin{equation}\label{uni1-1}
\|\UU^n\|_{1,r^-}^{r^-}+\|\SSS^n\|_{(r^+)'}^{(r^+)'}+\frac{1}{k}\|\UU^n\|^t_t\leq C_1,
\end{equation}
where $C_1$ is independent of $n$.

Next, we test with $C^n-c_d$ in \eqref{minfty3}, and by \eqref{Bfree} we obtain
\[\int_{\Omega}\q^n_c\cdot\nabla C^n\dx=\int_{\Omega}\q^n_c\cdot\nabla c_d\dx+B_c[C^n,\UU^n,c_d].\]
From \eqref{q1}, \eqref{q2}, H\"older's inequality and Young's inequality we have
\begin{align*}
\|\nabla C^n\|^2_2
&\leq C\int_{\Omega}|\nabla C^n||\nabla c_d|\dx+B_c[C^n,\UU^n,c_d]\\
&\leq\varepsilon\|\nabla C^n\|^2_2+C(\varepsilon)\|\nabla c_d\|^2_2+B_c[C^n,\UU^n,c_d].
\end{align*}
Furthermore, by Sobolev embedding,
\begin{align*}
B_c[C^{n},\UU^{n},c_d]
&=\frac{1}{2}\int_{\Omega}c_d\UU^{n}\cdot\nabla C^{n}\dx-\frac{1}{2}\int_{\Omega}C^{n}\UU^{n}\cdot\nabla c_d\dx\\
&=\int_{\Omega}c_d\UU^{n}\cdot\nabla C^{n}\dx+\frac{1}{2}\int_{\Omega}C^{n}({\rm{div}}\,\UU^{n})c_d\dx\\
&\leq\|c_d\|_{\infty}\|\UU^{n}\|_2\|\nabla C^{n}\|_2+\frac{\|c_d\|_{\infty}}{2}\|C^{n}\|_{(r^-)'}\|{\rm{div}}\,\UU^{n}\|_{r^-}\\
&\leq C\|\UU^{n}\|_{1,r^-}\|\nabla C^{n}\|_2+C\|\UU^{n}\|_{1,r^-}\|\nabla C^{n}\|_{\frac{3r^-}{4r^--3}}\\
&\leq C(\varepsilon)\|\UU^{n}\|^2_{1,r^-}+\varepsilon\|\nabla C^{n}\|^2_2.
\end{align*}
Hence, from \eqref{q1} and \eqref{uni1},
\begin{equation}\label{uni2}
\int_{\Omega}|\nabla C^n|^2+|\q_c^n|^2\dx\leq C_2,
\end{equation}
where $C_2$ is independent of $n$. Thus we have
\begin{equation}\label{uni2-1}
\|C^n\|^2_{1,2}+\|\q^n_c\|^2_2\leq C_2,
\end{equation}
where $C_2$ is independent of $n$.

Now, since $\frac{3}{2}>\frac{t}{t-2}$, by Sobolev embedding and the uniform estimates \eqref{uni1} and \eqref{uni2}, for $s>3$
sufficiently close to $3$,
\[\|C^n\UU^n\|_s\leq\|C^n\|_6\|\UU^n\|_{\frac{6s}{6-s}}\leq C\|C^n\|_{1,2}\|\UU^n\|_t\leq C,\]
where $C$ is independent of $n$.
Also, for $s>3$ sufficiently close to $3$, we have
\[\|\nabla C^n\cdot\UU^n\|_{\frac{3s}{s+3}}\leq\|\nabla C^n\|_2\|\UU^n\|_{\frac{6s}{6-s}}\leq C\|C^n\|_{1,2}\|\UU^n\|_t\leq C,\]
where $C$ is independent of $n$.

Therefore, we can apply Theorem \ref{DeGiorgi} with $\boldsymbol{F}=C^n\UU^n$ and $g=\nabla C^n\cdot\UU^n$. Hence, there exists
an $\alpha_1\in(0,1)$ such that
\begin{equation}\label{uni3}
\|C^n\|_{C^{0,\alpha_1}(\overline{\Omega})}\leq C_3.
\end{equation}
Since $C^{0,\alpha_1}(\overline{\Omega})\hookrightarrow\hookrightarrow C^{0,\tilde{\alpha}_1}(\overline{\Omega})$ for all $\tilde{\alpha}_1\in(0,\alpha_1)$, we have
\[C^n\rightarrow c\qquad{\rm{strongly}}\,\,\,{\rm{in}}\,\,\,C^{0,\tilde{\alpha}_1}(\overline{\Omega}),\]
which implies that
\[r\circ C^n\rightarrow r\circ c\qquad{\rm{strongly}}\,\,\,{\rm{in}}\,\,\,C^{0,\beta_1}(\overline{\Omega}),\]
for some $\beta_1\in(0,1)$.

We now apply Proposition \ref{disinfsup1}. For a given $r^+>0$, choose $j>\max\{r^+,2\}$. Then, since $r^->\frac{3}{2}$,
we have that $W^{1,j}_0(\Omega)^3\hookrightarrow L^{2(r^-)'}(\Omega)^3$ by Sobolev embedding.
Furthermore, since $\frac{t}{t-2}<r^-$, we have that $2(r^-)'<t$. Now, from \eqref{PGal1} and \eqref{B_u1},
\begin{align*}
\beta_r\|P^n\|_{j'}
&\leq\sup_{0\neq\VV\in\V^n}\frac{\langle{\rm{div}}\,\VV,P^n\rangle}{\|\VV\|_{1,j}}\\
&\leq\sup_{0\neq\VV\in\V^n}\frac{|\int_{\Omega}\SSS^n\cdot\DD\VV\dx+B_u[\UU^n,\UU^n,\VV]-\langle\boldsymbol{f},\VV\rangle|}{\|\VV\|_{1,j}}\\
&\leq C\sup_{0\neq\VV\in\V^n}\frac{\|\SSS^n\|_{(r^+)'}\|\VV\|_{1,r^+}}{\|\VV\|_{1,r^+}}\\
&\hspace{4mm}+C\sup_{0\neq\VV\in\V^n}\frac{{\color{black}{\|\UU^n\|^2_t}}\|\VV\|_{1,r^-}+\|\boldsymbol{f}\|_{-1}\|\VV\|_{1,r^-}}{\|\VV\|_{1,r^-}}\\
&\hspace{4mm}+C\sup_{0\neq\VV\in\V^n}\frac{\|\UU^n\|_{2(r-)'}\|\VV\|_{2(r-)'}\|\UU^n\|_{1,r^-}}{\|\VV\|_{2(r-)'}}.
\end{align*}
Hence, by noting \eqref{uni1},
\begin{equation}\label{uni4}
\|P^n\|_{j'}\leq C_4,
\end{equation}
where $C_4$ is independent of $n$.

Now, by \eqref{uni1}--\eqref{uni4}, thanks to the reflexivity of the relevant spaces
and by compact embedding, we can extract (not relabelled) subsequences such that
\begin{align}
\UU^n & \rightharpoonup \uu && {\rm{weakly}}\,\,{\rm{in}}\,\,W^{1,r^-}_{0}(\Omega)^3\cap L^t(\Omega)^3, \label{con1}\\
\UU^n & \rightarrow \uu && {\rm{strongly}}\,\,{\rm{in}}\,\,L^\sigma(\Omega)^3\qquad\forall\, \sigma\in{\color{black}{[1,t)}},\label{con2}\\
|\UU^n|^{t-2}\UU^n&\rightharpoonup|\uu|^{t-2}\uu&&{\rm{weakly}}\,\,{\rm{in}}\,\,L^{\frac{t}{t-1}}(\Omega)^3,\label{consub}\\
C^n & \rightharpoonup c && {\rm{weakly}}\,\,{\rm{in}}\,\,W^{1,2}(\Omega),\label{con3}\\
C^n & \rightarrow c && {\rm{strongly}}\,\,{\rm{in}}\,\,C^{0,\tilde{\alpha}_1}(\overline{\Omega}),\label{con4}\\
P^n & \rightharpoonup p && {\rm{weakly}}\,\,{\rm{in}}\,\,L^{j'}(\Omega)\qquad\forall\, j>\max\{r^+,2\},\label{con5}\\
\SSS^n & \rightharpoonup\bar{\SSS} && {\rm{weakly}}\,\,{\rm{in}}\,\,L^{(r^+)'}(\Omega)^{3\times 3},\label{con6}\\
\q_c^n & \rightharpoonup \bar{\q}_c && {\rm{weakly}}\,\,{\rm{in}}\,\,L^2(\Omega)^3.\label{con7}
\end{align}
Before proceeding further, we note that these limits, together with weak lower semicontinuity and  \eqref{uni1},
in conjunction with Korn's inequality, imply that
\begin{equation}\label{desired1}
\int_{\Omega}|\nabla \uu|^{r(c)}+|\bar{\SSS}|^{r'(c)}\dx\leq C,
\end{equation}
hence the limit function $\uu$ is, in fact, contained in the space $W^{1,r(c)}_0(\Omega)^3$; see \cite{KS2017} for the
details of the proof of this.

Next, we shall prove that the limit function $\uu$ is pointwise divergence-free. For an arbitrary $q\in C^{\infty}_0(\Omega)$,
by \eqref{minfty2},
\begin{align*}
0&=\int_{\Omega}(\Pi^n_{\Q}q)\,{\rm{div}}\,\UU^n\dx\\
&=\int_{\Omega}(\Pi^n_{\Q}q-q)\,{\rm{div}}\,\UU^n\dx+\int_{\Omega}q({\rm{div}}\,\UU^n-{\rm{div}}\,\uu)\dx+\int_{\Omega}q\,{\rm{div}}\,\uu\dx.
\end{align*}
The first term tends to zero by \eqref{uni1}, \eqref{q_conv} and the second term converges to zero by \eqref{con1}. Therefore,
\[\int_{\Omega}q\,{\rm{div}}\,\uu\dx=0\qquad{\rm{for}}\,\,\,{\rm{any}}\,\,\,q\in C^{\infty}_0(\Omega),\]
which implies that ${\rm{div}}\, \uu = 0$ a.e. on $\Omega$.

Now, we shall identify the limit of the convective term $B_u[\cdot,\cdot,\cdot]$ as follows. For an arbitrary $\vv\in W^{1,\infty}_0(\Omega)^3$, we define $\VV^n\defeq\Pi^n_{\rm{div}}\vv\in\V^n$. Then, by \eqref{v_conv}, we have
\begin{equation}\label{VV_n}
\VV^n\rightarrow \vv\qquad{\rm{strongly}}\,\,\,{\rm{in}}\,\,\,W^{1,\sigma}_0(\Omega)^2\,\,\,{\rm{for}}\,\,\,\sigma\in(1,\infty).
\end{equation}
By \eqref{con2},
\[\UU^n\otimes \UU^n\rightarrow \uu\otimes \uu\qquad{\rm{strongly}}\,\,\,{\rm{in}}\,\,\,L^{1+\varepsilon}(\Omega)^{3\times3}.\]
Hence, we can identify the second part of the convective term
\[-\int_{\Omega}(\UU^n\otimes \UU^n)\cdot\nabla \VV^n\dx\rightarrow-\int_{\Omega}(\uu\otimes \uu)\cdot\nabla \vv\dx\qquad{\rm{as}}\,\,\,n\rightarrow\infty.\]

Also, we assert that $\UU^n\cdot \VV^n\rightarrow \uu\cdot \vv$ strongly in $L^{(r^-)'}(\Omega)$. Indeed,
\begin{align*}
\|\UU^n\cdot \VV^n-\uu\cdot \vv\|_{(r^-)'}
&\leq\|(\VV^n-\vv)\UU^n+(\UU^n-\uu)\vv\|_{(r^-)'}\\
&\leq\|\VV^n-\vv\|_\sigma\|\UU^n\|_{t-\varepsilon}+\|\UU^n-\uu\|_{t-\varepsilon}\|\vv\|_\sigma
\end{align*}
for some $\sigma\in(1,\infty)$. The first term tends to zero thanks to \eqref{VV_n}, \eqref{con2} and the second term tends to zero by \eqref{con2}. Therefore, since ${\rm{div}}\,\uu=0$, we have
\begin{align*}
\int_{\Omega}(\UU^n\otimes \VV^n)\cdot\nabla \UU^n\dx
&=-\int_{\Omega}(\UU^n\otimes \UU^n)\cdot\nabla \VV^n\dx+\int_{\Omega}({\rm{div}}\,\UU^n)\,\UU^n\cdot \VV^n\dx\\
&\rightarrow-\int_{\Omega}(\uu\otimes \uu)\cdot\nabla \vv\dx\qquad{\rm{as}}\,\,\,n\rightarrow\infty.
\end{align*}
Altogether, we then deduce that
\begin{equation}\label{convectionconv}
\lim_{n\rightarrow\infty}B_u[\UU^n,\UU^n,\VV^n]=-\int_{\Omega}(\uu\otimes \uu)\cdot\nabla  \vv\dx.
\end{equation}

Now, we are ready to pass $n$ to infinity in the Navier--Stokes equations. Since $\Pi^n_{\rm{div}}$ is linear, by noting \eqref{minfty1}, we have
\begin{align*}
\langle{\rm{div}}\,\vv,P^n\rangle
&=\langle{\rm{div}}\,\VV^n,P^n\rangle+\langle{\rm{div}}\,(\vv-\VV^n),P^n\rangle\\
&=\int_{\Omega}\SSS(C^n,\DD\UU^n)\cdot \DD\VV^n+\frac{1}{k}|\UU^n|^{t-2}\UU^n\cdot\VV^n \dx-\langle \boldsymbol{f},\VV^n\rangle+B_u[\UU^n,\UU^n,\VV^n]\\
&\,\,\,\,\,\,+\langle{\rm{div}}\,(\vv-\VV^n),P^n\rangle\\
&\rightarrow\int_{\Omega}\bar{\SSS}\cdot \DD\vv+\frac{1}{k}|\uu|^{t-2}\uu\cdot\vv+{\rm{div}}(\uu\otimes \uu)\cdot \vv\dx-\langle \boldsymbol{f},\vv\rangle,
\end{align*}
where we have used \eqref{con5}, \eqref{con6}, \eqref{consub}, \eqref{VV_n} and \eqref{convectionconv}. Also, by \eqref{con5} again,
\[\langle{\rm{div}}\,\vv,P^n\rangle\rightarrow\langle{\rm{div}}\,\vv,p\rangle.\]
Collecting all the limits gives us
\begin{equation}\label{limiteq1}
\int_{\Omega}\bar{\SSS}\cdot \DD \vv+\frac{1}{k}|\uu|^{t-2}\uu\cdot\vv+{\rm{div}}\,(\uu\otimes \uu)\cdot \vv\dx-\langle{\rm{div}}\,\vv,p\rangle=\langle \boldsymbol{f},\vv\rangle\qquad\forall\, \vv\in W^{1,\infty}_0(\Omega)^3.
\end{equation}
With the same argument as above, we also have that
\begin{equation}\label{limiteq1-1}
\int_{\Omega}\bar{\SSS}\cdot \DD \vv+\frac{1}{k}|\uu|^{t-2}\uu\cdot\vv+{\rm{div}}\,(\uu\otimes \uu)\cdot \vv\dx=\langle \boldsymbol{f},\vv\rangle\qquad\forall\, \vv\in W^{1,\infty}_{0,{\rm{div}}}(\Omega)^3.
\end{equation}

Note that by Proposition \ref{contiinfsup2} and \eqref{desired1}, we have
\[p\in L^{r'(c)}_0(\Omega).\]

Now, let us investigate the limit of the convection-diffusion equation, \eqref{minfty3}. For an arbitrary but fixed $z\in W^{1,2}_0(\Omega)$,
we define $Z^n\defeq\Pi^n_{\Z}z\in\Z^n$. Thanks to \eqref{con2} and \eqref{con4},
\begin{align*}
\|C^n\UU^n-c\uu\|_2
&\leq\|(C^n-c)\UU^n\|_2+\|c(\UU^n-\uu)\|_2\\
&\leq\|C^n-c\|_{\infty}\|\UU^n\|_2+\|c\|_{\infty}\|\UU^n-\uu\|_2\rightarrow0.
\end{align*}
Moreover, by \eqref{con2}, \eqref{z_conv} and Sobolev embedding,
\begin{align*}
\|Z^n\UU^n-z\uu\|_2
&\leq\|(Z^n-z)\UU^n\|_2+\|z(\UU^n-\uu)\|_2\\
&\leq\|Z^n-z\|_6\|\UU^n\|_3
+\|z\|_6\|\UU^n-\uu\|_3\\
&\leq C\|Z^n-z\|_{1,2}\|\UU^n\|_3+C\|z\|_{1,2}\|\UU^n-\uu\|_3\rightarrow0.
\end{align*}
In other words,
\begin{align}
C^n\UU^n&\rightarrow c\uu\qquad{\rm{strongly}}\,\,\,{\rm{in}}\,\,\,L^2(\Omega)^3,\label{CU}\\
Z^n\UU^n&\rightarrow z\uu\qquad{\rm{strongly}}\,\,\,{\rm{in}}\,\,\,L^2(\Omega)^3.\label{ZU}
\end{align}
From \eqref{con3} and \eqref{ZU},
\begin{align*}
&\bigg|\int_{\Omega}Z^n\UU^n\cdot\nabla C^n\dx-\int_{\Omega}z\uu\cdot\nabla c\dx\bigg|\\
&\qquad\leq\int_{\Omega}|Z^n\UU^n-z\uu||\nabla C^n|\dx+\bigg|\int_{\Omega}z\uu\cdot(\nabla C^n-\nabla c)\dx\bigg|\rightarrow0.
\end{align*}
Therefore, as ${\rm{div}}\,\uu=0$ a.e. on $\Omega$, we obtain
\[\int_{\Omega}Z^n\UU^n\cdot\nabla C^n\dx\rightarrow\int_{\Omega}z\uu\cdot\nabla c\dx=-\int_{\Omega}c\uu\cdot\nabla z\dx\qquad{\rm{as}}\,\,\,n\rightarrow\infty.\]
Moreover, by \eqref{CU} and \eqref{z_conv},
\begin{align*}
&\,\,\,\,\,\,\,\,\bigg|\int_{\Omega}C^n\UU^n\cdot\nabla Z^n\dx-\int_{\Omega}c\uu\cdot\nabla z\dx\bigg|\\
&\qquad \leq\|C^n\UU^n\|_2\|Z^n-z\|_{1,2}+\|C^n\UU^n-c\uu\|_2\|z\|_{1,2}\rightarrow0.
\end{align*}
Altogether, we have
\[\lim_{n\rightarrow\infty} B_c[C^n,\UU^n,Z^n]=-\int_{\Omega}c\uu\cdot\nabla z\dx.\]
Finally, by \eqref{con7} and \eqref{z_conv}, we have
\[\int_{\Omega}\q_c(C^n,\nabla C^n,\DD\UU^n)\cdot\nabla Z^n\dx\rightarrow \int_{\Omega}\bar{\q}_c\cdot\nabla z\dx\qquad{\rm{as}}\,\,\,n\rightarrow\infty.\]
By collecting all the limits, we obtain that
\begin{equation}\label{limiteq2}
\int_{\Omega}\bar{\q}_c\cdot\nabla z-c\uu\cdot\nabla z\dx=0\qquad\forall\, z\in W^{1,2}_0(\Omega).
\end{equation}
As we can see from \eqref{limiteq1} and \eqref{limiteq2}, what we now need to prove is the identification of the limits:
\[\bar{\SSS}=\SSS(c,\DD\uu)\,\,\,\,\,{\rm{and}}\,\,\,\,\,\bar{\q}_c=\q_c(c,\nabla c, \DD\uu).\]
To this end, we require the following lemma.

{\color{black}{
\begin{lemma}\label{mainlem}
The sequences $\{\DD\UU^n\}_{n \in \mathbb{N}}$ and $\{C^n\}_{n \in \mathbb{N}}$ satisfy the following equality:
\begin{equation}\label{cptness}
\lim_{n\rightarrow\infty}\int_{\Omega}((\SSS(C^n,\DD\UU^n)-\SSS(C^n,\DD\uu))\cdot(\DD\UU^n-\DD\uu))^{\frac{1}{4}}\dx=0.
\end{equation}
\end{lemma}

The detailed proof of Lemma \ref{mainlem} is presented in Section 4.2 in \cite{KS2017}. Here, we shall briefly summarize the key steps of the proof as we shall require a similar, but more involved, argument in the next section. The strategy is to decompose the integral into several terms and to estimate them separately. To this end, for arbitrary but fixed $\chi>0$, we introduce the matrix-truncation function $T_{\chi}:\R^{3\times 3}\rightarrow\R^{3\times 3}$ by
\begin{displaymath}
T_{\chi}(\boldsymbol{M})=\left\{ \begin{array}{ll}
\boldsymbol{M} &\textrm{for $|\boldsymbol{M}|\leq\chi$,}\\
\chi\frac{\boldsymbol{M}}{|\boldsymbol{M}|}& \textrm{for $|\boldsymbol{M}|>\chi$.}
\end{array}\right.
\end{displaymath}

The essential step in the proof relies on using a discrete Lipschitz truncation technique. In  \cite{KS2017} a version of the discrete Lipschitz truncation method in variable-exponent norms was presented: see Theorem 3.15 in \cite{KS2017} and let $\UU^n_j$ denote the discrete Lipschitz truncation of the function of $\UU^n$.

The most important and difficult part of the proof is to estimate the following term:
\begin{equation}\label{diffesti}
\lim_{\chi \rightarrow 0}\lim_{j \rightarrow \infty}
\lim_{n \rightarrow \infty}\int_{\Omega}(\SSS(C^n,\DD\UU^n)-\SSS(C^n,T_{\chi}(\DD\uu)))\cdot(\DD\UU^n_j-T_{\chi}(\DD\uu))\dx\leq 0;
\end{equation}
(see eq. (4.23) in \cite{KS2017}).
The other terms arising from the decomposition can be easily estimated by using the uniform bound \eqref{uni1}, H\"older's inequality and the discrete Lipschitz truncation theorem, Theorem 3.15 in \cite{KS2017}.

To estimate \eqref{diffesti}, we introduce the following discretely divergence-free approximations with zero trace on $\partial\Omega$:
\begin{align*}
\boldsymbol{\Psi}^n_j&\defeq\mathcal{B}^n({\rm{div}}\,\UU^n_j),\\
\boldsymbol{\Phi}^n_j&\defeq \UU^n_j-\boldsymbol{\Psi}^n_j.
\end{align*}
Here $\mathcal{B}^n$ is a discrete Bogovski\u{\i} operator defined in Section 3.4 of \cite{KS2017}. It is then clear that $\boldsymbol{\Phi}^n_j$ has zero trace on $\partial\Omega$ and, by construction, $\boldsymbol{\Phi}^n_j\in\V^n_{\rm{div}}$. Moreover, it can be easily verified, by using basic properties of the discrete Lipschitz truncation and the discrete Bogovski\u{\i} operator, that
\begin{align}
\boldsymbol{\Phi}^n_j&\rightharpoonup \UU_j-\mathcal{B}({\rm{div}}\,\UU_j)\eqdef\boldsymbol{\Phi}_j&&{\rm{weakly}}\,\,{\rm{in}}\,\, W^{1,\sigma}_0(\Omega)^3,\label{main5.171}\\
\boldsymbol{\Phi}^n_j&\rightarrow\boldsymbol{\Phi}_j&&{\rm{strongly}}\,\,{\rm{in}}\,\,L^{\sigma}(\Omega)^3,\label{main5.181}
\end{align}
as $n\rightarrow\infty$, where $\sigma\in(1,\infty)$ is arbitrary. We can then rewrite \eqref{diffesti} above in terms of this approximation to obtain
\begin{align*}
&\int_{\Omega}(\SSS(C^n,\DD\UU^n)-\SSS(C^n,T_{\chi}(\DD\uu)))\cdot(\DD\UU^n_j-T_{\chi}(\DD\uu))\dx\\
&=\int_{\Omega}\SSS(C^n,\DD\UU^n)\cdot(\DD\boldsymbol{\Phi}^n_j+\DD\boldsymbol{\Psi}^n_j)\dx\\
&\,\,\,\,\,\,-\int_{\Omega}\SSS(C^n,\DD\UU^n)\cdot T_{\chi}(\DD\uu)\dx-\int_{\Omega}\SSS(C^n,T_{\chi}(\DD\uu))\cdot(\DD\UU^n_j-T_{\chi}(\DD\uu))\dx\\
&\eqdef B^{n,1}_{\chi,j}-B^{n,2}_{\chi,j}-B^{n,3}_{\chi,j}.
\end{align*}
Now we use \eqref{minfty1} with $\VV=\boldsymbol{\Phi}^n_j\in\V^n_{\rm{div}}$ and pass to the limit; thus we have, by \eqref{limiteq1-1}, that
\begin{align}
\lim_{n\rightarrow\infty}\int_{\Omega}\SSS^n\cdot \DD\boldsymbol{\Phi}^n_j\dx
&=-\lim_{n\rightarrow\infty}B_u[\UU^n,\UU^n,\boldsymbol{\Phi}^n_j]-\int_{\Omega}\frac{1}{k}|\UU^n|^{t-2}\UU^n\cdot\boldsymbol{\Phi}^n_j\dx+\lim_{n\rightarrow\infty}\langle \boldsymbol{f},\boldsymbol{\Phi}^n_j\rangle\label{important1}\\
&=\int_{\Omega}(\uu\otimes \uu)\cdot\nabla\boldsymbol{\Phi_j}-\frac{1}{k}|\uu|^{t-2}\uu\cdot\boldsymbol{\Phi_j}\dx+\langle \boldsymbol{f},\boldsymbol{\Phi_j}\rangle\label{important2}\\
&=\int_{\Omega}\bar{\SSS}\cdot \DD\boldsymbol{\Phi_j}\dx.\label{important3}
\end{align}

Furthermore, with the help of Lipschitz truncation, we can show that
\begin{align}
\lim_{n\rightarrow\infty}\int_{\Omega}\SSS^n\cdot\DD\boldsymbol{\Psi}^n_j\dx&\leq\left(\frac{C}{2^{j/r^+}}\right)^{\gamma(r^-,r^+)},\label{important4}\\
\int_{\Omega}\bar{\SSS}\cdot\DD\mathcal{B}({\rm{div}}\,\UU_j)\dx&\leq\left(\frac{C}{2^{j/r^+}}\right)^{\gamma(r^-,r^+)}.\label{important5}
\end{align}
Altogether, we have
\[\lim_{\chi\rightarrow\infty}\lim_{j\rightarrow\infty}\lim_{n\rightarrow\infty}\left(B^{n,1}_{\chi,j}-B^{n,2}_{\chi,j}-B^{n,3}_{\chi,j}\right)\leq\lim_{\chi\rightarrow\infty}\int_{\Omega}(\bar{\SSS}-\SSS(c,T_{\chi}(\DD\uu)))\cdot(\DD\UU_j-T_{\chi}(\DD\uu))\dx\]
The last limit is equal to zero by using the Dominated Convergence Theorem. That completes the proof of \eqref{diffesti}, and thereby also of the most technical step in the proof of the lemma.}}

\medskip

 Now we are ready to identify the limits. In the above lemma, since the integrand is nonnegative, \eqref{cptness} also holds with $\Omega$ replaced by the set
$Q_{\gamma}\subset\Omega$ defined by
\[Q_{\gamma}\defeq\{x\in\Omega:|\DD\uu|\leq\gamma\},\]
with a given  $\gamma>0$; thus, from the sequence of integrands featuring in  \eqref{cptness}, we can extract a subsequence
(again not relabelled), which converges to zero almost everywhere in $Q_{\gamma}$. Then, by Egoroff's Theorem,
for an arbitrary $\varepsilon>0$, there exists a subset $Q^{\varepsilon}_{\gamma}\subset Q_{\gamma}\subset \Omega$ satisfying
$|Q_{\gamma}\setminus Q^{\varepsilon}_{\gamma}|<\varepsilon$, where the convergence of  integrands is uniform.
Note that, thanks to the choice of $Q^{\varepsilon}_{\gamma}$, we have
\[\lim_{\gamma\rightarrow\infty}\lim_{\varepsilon\rightarrow0}|\Omega\setminus Q^{\varepsilon}_{\gamma}|=
\lim_{\gamma\rightarrow\infty}\lim_{\varepsilon\rightarrow0}\left[ |\Omega\setminus Q_{\gamma}| + |Q_{\gamma}\setminus Q^{\varepsilon}_{\gamma}|
\right]= 0.\]
Moreover, we have from the uniform convergence of the integrands that
\begin{equation}\label{main5.21}
\lim_{n\rightarrow\infty}\int_{Q^{\varepsilon}_{\gamma}}(\SSS(C^n,\DD\UU^n)-\SSS(C^n,\DD\uu))\cdot (\DD\UU^n-\DD\uu)\dx=0.
\end{equation}
Since $\DD\uu$ is bounded on $Q^{\varepsilon}_{\gamma}$, by the Dominated Convergence Theorem we have $\SSS(C^n,\DD\uu)\rightarrow \SSS(c,\DD\uu)$ strongly in $L^q(\Omega)^{3\times 3}$ for any $q\in[1,\infty)$.
Hence, from the above $L^q$-convergence, \eqref{conv1}, and \eqref{main5.21}, we obtain
\[\lim_{n\rightarrow\infty}\int_{Q^{\varepsilon}_{\gamma}}\SSS(C^n,\DD\UU^n)\cdot(\DD\UU^n-\DD\uu)\dx=0.\]
Thus, by the boundedness of $\DD\uu$ on $Q^{\varepsilon}_{\gamma}$ and  \eqref{con6}, we have
\begin{equation}\label{main5.22}
\lim_{n\rightarrow\infty}\int_{Q^{\varepsilon}_{\gamma}}\SSS(C^n,\DD\UU^n)\cdot \DD\UU^n\dx=\int_{Q^{\varepsilon}_{\gamma}}\bar{\SSS}\cdot \DD\uu\dx.
\end{equation}
Now, let $\boldsymbol{B}\in L^{\infty}(Q^{\varepsilon}_{\gamma})^{3\times 3}$ be arbitrary but fixed.
From the monotonicity \eqref{S2}, \eqref{main5.22}, the $L^q$-convergence of $\SSS(C^n,\boldsymbol{B})\rightarrow \SSS(c,\boldsymbol{B})$ and the weak convergence \eqref{conv1}, we have
\begin{align*}
0&\leq\lim_{n\rightarrow\infty}\int_{Q^{\varepsilon}_{\gamma}}(\SSS(C^n,\DD\UU^n)-\SSS(C^n,\boldsymbol{B}))\cdot(\DD\UU^n-\boldsymbol{B})\dx\\
&=\int_{Q^{\varepsilon}_{\gamma}}\bar{\SSS}\cdot(\DD\uu-\boldsymbol{B})\dx-\int_{Q^{\varepsilon}_{\gamma}}\SSS(c,\boldsymbol{B})\cdot(\DD\uu-\boldsymbol{B})\dx\\
&=\int_{Q^{\varepsilon}_{\gamma}}(\bar{\SSS}-\SSS(c,\boldsymbol{B}))\cdot(\DD\uu-\boldsymbol{B})\dx.
\end{align*}
Now we are ready to use Minty's trick. First, we choose $\boldsymbol{B}=\DD\uu\pm\lambda \boldsymbol{A}$ with $\lambda>0$ and $\boldsymbol{A}\in L^{\infty}(Q^{\varepsilon}_{\gamma})^{3\times 3}$. Then, passing to the limit $\lambda\rightarrow0$,
the continuity of $\SSS$ gives us
\[\int_{Q^{\varepsilon}_{\gamma}}(\bar{\SSS}-\SSS(c,\DD\uu))\cdot \boldsymbol{A}\dx=0.\]
Hence, we have that
\[\bar{\SSS}=\SSS(c,\DD\uu)\,\,\,{\rm{a.e.}}\,\,\,{\rm{on}}\,\,\,Q^{\varepsilon}_{\gamma}.\]
Now we pass $\varepsilon\rightarrow0$ and then $\gamma\rightarrow\infty$ to conclude that
\begin{equation}\label{ident}
\bar{\SSS}=\SSS(c,\DD\uu)\,\,\,{\rm{a.e.}}\,\,\,{\rm{on}}\,\,\,\Omega.
\end{equation}

Finally, since $\SSS$ is strictly monotonic and $C^n\rightarrow c$ in $C^{0,\tilde{\alpha}_1}(\overline{\Omega})$,
by \eqref{cptness} we deduce that
\begin{equation}\label{a.e.conv}
\DD\UU^n\rightarrow \DD\uu\,\,\,{\rm{a.e.}}\,\,\,{\rm{on}}\,\,\,\Omega.
\end{equation}
By the Dominated Convergence Theorem, with \eqref{con3}, \eqref{con4} and \eqref{a.e.conv}, we obtain that
\[\q_c(C^n,\nabla C^n,\DD\UU^n)\rightharpoonup \q_c(c,\nabla c, \DD\uu)\qquad{\rm{weakly}}\,\,{\rm{in}}\,\,L^2(\Omega)^3.\]
Therefore, by the uniqueness of the weak limit, we can identify
\begin{equation}\label{ident2}
\bar{\q}_c=\q_c(c,\nabla c,\DD\uu).
\end{equation}
\end{subsection}
\end{section}

\begin{section}{Proof of Theorem \ref{mainthm2}}
\begin{subsection}{Minimum and maximum principles}
Before we proceed, let us prove minimum and maximum principles for the concentration. Let $\varphi^k_1=(c^k-\min_{x\in\partial\Omega}c_d)_-$
and $\varphi^k_2=(c^k-\max_{x\in\partial\Omega}c_d)_+$. Since $c^k=c_d$ on $\partial\Omega$, it is clear that
$\varphi^k_1,\varphi^k_2\in W^{1,2}_0(\Omega)$, so we can test with $\varphi^k_1$ and $\varphi^k_2$ in \eqref{limiteq2}. Therefore, we have
\begin{equation}\label{min}
-\int_{\Omega}\uu^kc^k\cdot\nabla\varphi^k_1\dx+\int_{\Omega}\bar{\q}_c\nabla\varphi^k_1\dx=0,
\end{equation}
\begin{equation}\label{max}
-\int_{\Omega}\uu^kc^k\cdot\nabla\varphi^k_2\dx+\int_{\Omega}\bar{\q}_c\nabla\varphi^k_2\dx=0.
\end{equation}
We first consider \eqref{min}. From \eqref{q2} with integration by parts we obtain
\[\int_{\Omega^-}\uu^k\cdot\nabla c^k\varphi^k_1\dx+\int_{\Omega^-}C|\nabla c^k|^2\dx\leq0,\]
where $\Omega^-=\{x\in\Omega:\varphi^k_1(x)<0\},$ since ${\rm{div}}\,\uu^k=0$ and $\uu^k=0$ on $\partial\Omega$. By using the fact that $\nabla c^k=\nabla\varphi^k_1$ on $\Omega^-$ and the extension of $\nabla c^k$ from $\Omega^-$ to the whole domain $\Omega$ by using the negative part, we have
\[\int_{\Omega}\uu^k\cdot\nabla\varphi^k_1 \varphi^k_1\dx+\int_{\Omega}C|\nabla\varphi^k_1|^2\dx\leq0.\]
Note that
\[\int_{\Omega}\uu^k\cdot\nabla\varphi^k_1\varphi^k_1\dx=\frac{1}{2}\int_{\Omega}\uu^k\cdot\nabla|\varphi^k_1|^2\dx=-\frac{1}{2}\int_{\Omega}({\rm{div}}\,\uu^k)|\varphi^k_1|^2\dx=0,\]
and thus,
\[\varphi^k_1=(c^k-\min_{x\in\partial\Omega}c_d)_-={\rm{constant}}\,\,{\rm{a.e.}}\,\,{\rm{in}}\,\,\Omega.\]
In the same way, we can also show that
\[\varphi^k_2=(c^k-\max_{x\in\partial\Omega}c_d)_+={\rm{constant}}\,\,{\rm{a.e.}}\,\,{\rm{in}}\,\,\Omega.\]
By combining the above results we finally obtain that
\begin{equation}\label{min_max}
\min_{x\in\partial\Omega}c_d\leq c^k\leq\max_{x\in\partial\Omega}c_d\qquad{\rm{a.e.}}\,\,\,{\rm{in}}\,\,\,\Omega.
\end{equation}
\end{subsection}

\begin{subsection}{The limit $k\rightarrow\infty$}
First, note that by weak lower semicontinuity of the norm-function, and \eqref{uni1-1}, \eqref{uni2-1} and \eqref{uni4}, we obtain the following
uniform estimates, independent of $k\in\mathbb{N}$:
\begin{equation}\label{kuni1}
\|\uu^k\|^{r^-}_{1,r^-}+\|\SSS(c^k,\DD\uu^k)\|^{(r^+)'}_{(r^+)'}+\frac{1}{k}\|\uu^k\|^t_t\leq C_1,
\end{equation}
\begin{equation}\label{kuni2}
\|c^k\|^2_{1,2}+\|\q_c(c^k,\nabla c^k,\DD\uu^k)\|^2_2\leq C_2,
\end{equation}
\begin{equation}\label{kuni3}
\|p^k\|^{j'}_{j'}\leq C_3,
\end{equation}
for some positive constants $C_1$, $C_2$ and $C_3$, which are independent of $k\in\mathbb{N}$.

Now, since $r^->\frac{3}{2}$, by the min/max principle \eqref{min_max}, Sobolev embedding and the uniform estimate \eqref{kuni1},
for $s>3$ sufficiently close to $3$,
\[\|c^k\uu^k\|_s\leq\|c^k\|_{\infty}\|\uu^k\|_s\leq C\|\uu^k\|_{1,r^-}\leq C.\]

Therefore, we can again apply Theorem \ref{DeGiorgi} with $\boldsymbol{F}=c^k\uu^k$ and $g=0$. Hence, there exists an $\alpha_2\in(0,1)$ such that
\begin{equation}\label{kuni4}
\|c^k\|_{C^{0,\alpha_2}(\overline{\Omega})}\leq C_4,
\end{equation}
for some positive constant $C_4$ independent of $k\in\mathbb{N}$. Since $C^{0,\alpha_2}(\overline{\Omega})\hookrightarrow\hookrightarrow C^{0,\tilde{\alpha}_2}(\overline{\Omega})$ for all $\tilde{\alpha}_2\in(0,\alpha_2)$, we have
\[c^k\rightarrow c\qquad{\rm{strongly}}\,\,\,{\rm{in}}\,\,\,C^{0,\tilde{\alpha}_2}(\overline{\Omega}),\]
which implies that
\[r\circ c^k\rightarrow r\circ c\qquad{\rm{strongly}}\,\,\,{\rm{in}}\,\,\,C^{0,\beta_2}(\overline{\Omega}),\]
for some $\beta_2\in(0,1)$.

Therefore, by the reflexivity of the relevant spaces and compact embedding, there exists a subsequence (not relabelled) such that
\begin{align}
\uu^k & \rightharpoonup \uu && {\rm{weakly}}\,\,{\rm{in}}\,\,W^{1,r^-}_{0,{\rm{div}}}(\Omega)^3, \label{co1}\\
\uu^k & \rightarrow \uu && {\rm{strongly}}\,\,{\rm{in}}\,\,L^{2(1+\varepsilon)}(\Omega)^3,\label{co2}\\
c^k & \rightharpoonup c && {\rm{weakly}}\,\,{\rm{in}}\,\,W^{1,2}(\Omega),\label{co3}\\
c^k & \rightarrow c && {\rm{strongly}}\,\,{\rm{in}}\,\,C^{0,\tilde{\alpha}_2}(\overline{\Omega}),\label{co4}\\
p^k & \rightharpoonup p && {\rm{weakly}}\,\,{\rm{in}}\,\,L^{j'}(\Omega)\qquad\forall\, j>\max\{r^+,2\},\label{co5}\\
\SSS(c^k,\DD\uu^k) & \rightharpoonup\hat{\SSS} && {\rm{weakly}}\,\,{\rm{in}}\,\,L^{(r^+)'}(\Omega)^{3\times 3},\label{co6}\\
\q_c(c^k,\nabla c^k, \DD\uu^k) & \rightharpoonup \hat{\q}_c && {\rm{weakly}}\,\,{\rm{in}}\,\,L^2(\Omega)^3.\label{co7}
\end{align}

Again, by the weak lower semicontinuity of norms, \eqref{desired1} and \eqref{co4} together with Korn's inequality, we have that
\begin{equation}\label{desired2}
\int_{\Omega}|\nabla \uu|^{r(c)}+|\hat{\SSS}|^{r'(c)}\dx\leq C,
\end{equation}
and thus the weak solution $\uu$ is in the desired space $W^{1,r(c)}_0(\Omega)^3$.

Now we shall let $k\rightarrow\infty$ in \eqref{limiteq1}, with $\vv\in W^{1,\infty}_0(\Omega)^3$ chosen arbitrarily. By \eqref{co2},
\[\uu^k\otimes \uu^k\rightarrow \uu\otimes \uu\qquad{\rm{strongly}}\,\,\,{\rm{in}}\,\,\,L^{1+\varepsilon}(\Omega)^{3\times3}.\]
Thus, we can identify the limit of the convective term
\[-\int_{\Omega}(\uu^k\otimes \uu^k)\cdot\nabla \vv\dx\rightarrow-\int_{\Omega}(\uu\otimes \uu)\cdot\nabla \vv\dx\qquad{\rm{as}}
\,\,\,k\rightarrow\infty,\qquad \forall\,\vv\in W^{1,\infty}_0(\Omega)^3.\]

Next, by \eqref{kuni1}, we have that
\[\frac{1}{k}\|\uu^k\|^{t-1}_t\rightarrow0\qquad{\rm{as}}\,\,\,k\rightarrow\infty.\]
Therefore, we have
\[\frac{1}{k}\bigg|\int_{\Omega}|\uu^k|^{t-2}\uu^k\cdot\vv\dx\bigg|\leq\frac{1}{k}\|\uu^k\|^{t-1}_t\|\vv\|_t
\rightarrow0\qquad{\rm{as}}\,\,\,k\rightarrow\infty,\qquad \forall\,\vv\in W^{1,\infty}_0(\Omega)^3.\]

We recall from the identification asserted in \eqref{ident} that $\bar{\SSS}=\SSS(c,\DD\uu)$ a.e. on $\Omega$; more precisely, with the index $k$ reinstated in our notation, $\bar{\SSS}^k=\SSS(c^k,\DD\uu^k)$ a.e. on $\Omega$. Hence, from \eqref{co6} and \eqref{co5}, we obtain
\[
\langle{\rm{div}}\,\vv,p^k\rangle\rightarrow\langle{\rm{div}}\,\vv,p\rangle\quad{\rm{and}}\quad
\int_{\Omega}\bar{\SSS}^k\cdot\DD\vv\dx\rightarrow\int_{\Omega}\hat{\SSS}\cdot\DD\vv\dx\quad{\rm{as}}\,\,\,k\rightarrow\infty,\quad \forall\,\vv\in W^{1,\infty}_0(\Omega)^3.
\]

Altogether, we have
\begin{equation}\label{fianl1}
\int_{\Omega}\hat{\SSS}\cdot \DD \vv+(\uu\otimes \uu)\cdot \nabla\vv\dx-\langle{\rm{div}}\,\vv,p\rangle=\langle \boldsymbol{f},\vv\rangle\qquad\forall\, \vv\in W^{1,\infty}_0(\Omega)^3.
\end{equation}
Furthermore, it is clear that
\begin{equation}\label{final1-1}
\int_{\Omega}\hat{\SSS}\cdot \DD \vv+(\uu\otimes \uu)\cdot \nabla\vv\dx=\langle \boldsymbol{f},\vv\rangle\qquad\forall\, \vv\in W^{1,\infty}_{0,{\rm{div}}}(\Omega)^3.
\end{equation}
Note that by Proposition \ref{contiinfsup2} and \eqref{desired2} we have
\[p\in L^{r'(c)}_0(\Omega).\]

Now, let us investigate the limit of the concentration equation \eqref{limiteq2}.
Let us choose an arbitrary, but fixed, $z\in W^{1,2}_0(\Omega)$. By \eqref{co2} and \eqref{co4},
\[\|c^k\uu^k-c\uu\|_2
\leq\|(c^k-c)\uu^k\|_2+\|c(\uu^k-\uu)\|_2
\leq\|c^k-c\|_{\infty}\|\uu^k\|_2+\|c\|_{\infty}\|\uu^k-\uu\|_2\rightarrow0.\]
In other words,
\[c^k\uu^k\rightarrow c\uu\qquad{\rm{strongly}}\,\,\,{\rm{in}}\,\,\,L^2(\Omega)^3.\label{c_u}\]
Hence we have
\[\int_{\Omega}c^k\uu^k\cdot\nabla z\dx\rightarrow\int_{\Omega}c\uu\cdot\nabla z\dx.\]
Recalling the identification \eqref{ident2} and reinstating the index $k$, we have
$\bar{\q}_c^k:=\q_c(c^k,\nabla c^k, \DD\uu^k)$; hence, by \eqref{co7}, we get
\[\int_{\Omega}\bar{\q}_c^k\cdot\nabla z\dx\rightarrow \int_{\Omega}\hat{\q}_c\cdot\nabla z\dx\qquad{\rm{as}}\,\,\,k\rightarrow\infty.\]
By collecting the above limits, we deduce that
\begin{equation}\label{final2}
\int_{\Omega}\hat{\q}_c\cdot\nabla z-c\uu\cdot\nabla z\dx=0\qquad\forall\, z\in W^{1,2}_0(\Omega).
\end{equation}

As a final step, we need to identify the limits:
\[\hat{\SSS}=\SSS(c,\DD\uu)\,\,\,\,\,{\rm{and}}\,\,\,\,\,\hat{\q}_c=\q_c(c,\nabla c, \DD\uu).\]
{\color{black}{To this end, analogously as before, we need to prove the following equality:
\begin{equation}\label{cptness_k}
\lim_{k\rightarrow\infty}\int_{\Omega}((\SSS(c^k,\DD\uu^k)-\SSS(c^k,\DD\uu))\cdot(\DD\uu^k-\DD\uu))^{\frac{1}{4}}\dx=0.
\end{equation}
The proof is similar to the one presented in the previous section. The only part of the argument that we shall give here in detail is the proof of the analogue of \eqref{important1}--\eqref{important3} since we now have a different weak formulation at this level. The other parts of the proof proceed as Section 4.2 in \cite{KS2017}.

First we define a divergence-free approximation with zero trace as follows:
\[\boldsymbol{\Phi}^k_j\defeq\uu^k_j-\mathcal{B}({\rm{div}}\,\uu^k_j),\]
where $\mathcal{B}$ is the Bogovski\u{\i} operator introduced in Theorem \ref{contibog}. Then, as before, we have
\begin{align}
\boldsymbol{\Phi}^k_j&\rightharpoonup \uu_j-\mathcal{B}({\rm{div}}\,\uu_j)\eqdef\boldsymbol{\Phi}_j&&{\rm{weakly}}\,\,{\rm{in}}\,\, W^{1,\sigma}_0(\Omega)^3,\label{main5.1711}\\
\boldsymbol{\Phi}^k_j&\rightarrow\boldsymbol{\Phi}_j&&{\rm{strongly}}\,\,{\rm{in}}\,\,L^{\sigma}(\Omega)^3,\label{main5.1811}
\end{align}
as $k\rightarrow\infty$, where $\sigma\in(1,\infty)$ is arbitrary.

Let us further define $\boldsymbol{\chi}^{n,k}_{1,j}\defeq\Pi^n_{\rm{div}}\boldsymbol{\Phi}^k_j$. Then, by \eqref{v_conv},
\[\boldsymbol{\chi}^{n,k}_{1,j}\rightarrow\boldsymbol{\Phi}^k_j\qquad{\rm{strongly}}\,\,\,{\rm{in}}\,\,\,W^{1,\sigma}_0(\Omega)^3,\,\,\,\forall\,\sigma\in(1,\infty).\]
Now, by \eqref{minfty1},
\[\int_{\Omega}\SSS^n\cdot\DD\boldsymbol{\chi}^{n,k}_{1,j}\dx=-B_u[\UU^n,\UU^n,\boldsymbol{\chi}^{n,k}_{1,j}]-\int_{\Omega}\frac{1}{k}|\UU^n|^{t-2}\UU^n\cdot\boldsymbol{\chi}^{n,k}_{1,j}\dx+\langle\boldsymbol{f},\boldsymbol{\chi}^{n,k}_{1,j}\rangle.\]
If we take $n\rightarrow\infty$ in the above equality, we have
\begin{equation}\label{import6}
\int_{\Omega}\SSS(c^k,\DD\uu^k)\cdot\DD\boldsymbol{\Phi}^k_j\dx=\int_{\Omega}(\uu^k\otimes\uu^k)\cdot\nabla\boldsymbol{\Phi}^k_j-\frac{1}{k}|\uu^k|^{t-2}\uu^k\cdot\boldsymbol{\Phi}^k_j\dx+\langle\boldsymbol{f},\boldsymbol{\Phi}^k_j\rangle.
\end{equation}

Next, we define $\boldsymbol{\chi}^{n,k}_{2,j}\defeq\Pi^n_{\rm{div}}\boldsymbol{\Phi}_j$, and then we have
\[\boldsymbol{\chi}^{n,k}_{2,j}\rightarrow\boldsymbol{\Phi}_j\qquad{\rm{strongly}}\,\,\,{\rm{in}}\,\,\,W^{1,\sigma}_0(\Omega)^3,\,\,\,\forall\,\sigma\in(1,\infty).\]
Again, by \eqref{minfty1},
\[\int_{\Omega}\SSS^n\cdot\DD\boldsymbol{\chi}^{n,k}_{2,j}\dx=-B_u[\UU^n,\UU^n,\boldsymbol{\chi}^{n,k}_{2,j}]-\int_{\Omega}\frac{1}{k}|\UU^n|^{t-2}\UU^n\cdot\boldsymbol{\chi}^{n,k}_{2,j}\dx+\langle\boldsymbol{f},\boldsymbol{\chi}^{n,k}_{2,j}\rangle.\]
If we take $n\rightarrow\infty$, we have
\[\int_{\Omega}\SSS(c^k,\DD\uu^k)\cdot\DD\boldsymbol{\Phi}_j\dx=\int_{\Omega}(\uu^k\otimes\uu^k)\cdot\nabla\boldsymbol{\Phi}_j-\frac{1}{k}|\uu^k|^{t-2}\uu^k\cdot\boldsymbol{\Phi}_j\dx+\langle\boldsymbol{f},\boldsymbol{\Phi}_j\rangle.\]
Subsequently, if we pass $k$ to the infinity, we obtain
\begin{equation}\label{import7}
\int_{\Omega}\hat{\SSS}\cdot\DD\dx=\int_{\Omega}(\uu\otimes\uu)\cdot\nabla\boldsymbol{\Phi}_j+\langle\boldsymbol{f},\boldsymbol{\Phi}_j\rangle.
\end{equation}

Therefore, from \eqref{import6} and \eqref{import7}, we deduce that
\begin{align*}
\lim_{k\rightarrow\infty}\int_{\Omega}\SSS(c^k,\DD\uu^k)\cdot\DD\boldsymbol{\Phi}^k_j\dx
&=\lim_{k\rightarrow\infty}\int_{\Omega}(\uu^k\otimes\uu^k)\cdot\nabla\boldsymbol{\Phi}^k_j-\frac{1}{k}|\uu^k|^{t-2}\uu^k\cdot\boldsymbol{\Phi}^k_j\dx+\lim_{k\rightarrow\infty}\langle\boldsymbol{f},\boldsymbol{\Phi}^k_j\rangle\\
&=\int_{\Omega}(\uu\otimes\uu)\cdot\nabla\boldsymbol{\Phi}_j\dx+\langle\boldsymbol{f},\boldsymbol{\Phi}_j\rangle\\
&=\int_{\Omega}\hat{\SSS}\cdot\DD\boldsymbol{\Phi}_j\dx,
\end{align*}
which is the desired analogue of \eqref{important1}--\eqref{important3} corresponding to the limit $k\rightarrow \infty$, and thereby the proof of \eqref{cptness_k} has been completed.}}

We can then use the same argument as the one we employed in the previous section to
identify $\bar{\SSS} = \bar{\SSS}^k=\SSS(c^k,\DD\uu^k)$ and $\bar{\q}_c=\bar{\q}_c^k = \q_c(c^k,\nabla c^k,\DD\uu^k)$
(cf. \eqref{ident} and \eqref{ident2}, with the index $k$ reinstated),
and thus we can again identify $\hat{\SSS}=\SSS(c,\DD\uu)$, $\hat{\q}_c=\q_c(c,\nabla c,\DD\uu)$.
That completes the proof of the convergence theorem.
\end{subsection}
\end{section}

\section{Conclusions}

We have considered a system of nonlinear partial differential equations modelling the motion of an incompressible chemically reacting generalized Newtonian fluid in three space dimensions. The governing system consists of a steady convection-diffusion equation for the concentration and a generalized steady power-law-type fluid flow model for the velocity and the pressure,  where the viscosity depends on both the shear-rate and the concentration through a concentration-dependent power-law index.
We performed a rigorous convergence analysis of a finite element approximation of a regularized counterpart of the model;
specifically, we showed the convergence of the finite element method to a weak solution of the regularized model. We then
proved that weak solutions of the regularized problem converge to a weak solution of the original problem.

\section*{Acknowledgements}
Seungchan Ko's work was supported by the UK Engineering and Physical Sciences Research Council [EP/L015811/1].

\bibliography{references}
\bibliographystyle{abbrv}


\end{document}